\documentclass{amsart}

\usepackage{amssymb}
\usepackage{xcolor}
\usepackage{amsmath}
\usepackage{bm}

\newcommand{\half}{\frac{1}{2}}
\newcommand{\thalf}{\tfrac{1}{2}}
\newcommand{\summ}{\mathop{{\sum}^{\star}}}

\newcommand{\sym}{{\rm sym}}
\newcommand{\intt}{\int_{-\infty}^{\infty}}

\numberwithin{equation}{section}

\newtheorem{theorem}{Theorem}[section]

\newtheorem{lemma}[theorem]{Lemma}

\makeatletter
\@namedef{subjclassname@2010}{%
  \textup{2010} Mathematics Subject Classification}
\makeatother

\begin{document}

\title{Subconvexity bounds for twisted $L$-functions, II}

\author{Rizwanur Khan}
\address{
Department of Mathematics\\ University of Mississippi\\ University, MS 38677}
\email{rrkhan@olemiss.edu }

\thanks{The author was supported by the National Science Foundation Grant DMS-2001183 and the Simons Foundation (award 630985). Any opinions, findings and conclusions or recommendations expressed in this material are those of the authors and do not necessarily reflect the views of the National Science Foundation.}

\begin{abstract} 
We prove hybrid subconvexity bounds twisted $L$-functions $L(s,f\times \chi)$ at the central point using a fourth moment estimate, including a new instance of the Burgess subconvexity bound.
\end{abstract}

\maketitle

\section{Introduction}

Let $f$ be a holomorphic Hecke cusp form for the congruence group $\Gamma_0(q)$, and $\chi$ a primitive Dirichlet character of modulus $p$, where $q$ and $p$ are distinct odd primes. The twisted $L$-function $L(s, f\times \chi)$ corresponds by Atkin-Lehner theory to a newform of level $qp^2$ and nebentypus $\chi^2$. Thus the analtyic conductor of this $L$-function in the level aspect is $qp^2$, and the (hybrid) subconvexity problem is to obtain the bound
\[
L(\thalf, f\times \chi)\ll (qp^2)^{\kappa+\epsilon}
\]
for some $0\le \kappa< \frac14$, where we have normalized so that $s=\half$ is the central point. Such a bound is known for all ranges of $q$ and $p$, by work of Blomer, Harcos, and Michel \cite{bloharmic2}. Their subconvexity bound ($\kappa=\frac14-\frac{1}{1889}$) has the advantage of being very general, but much stronger bounds are known in certain ranges. Such bounds are contained in various works, each of which has its unique strength. Blomer and Harcos \cite{blohar} proved the subconvexity bound 
\begin{align}
\label{blo-har} L(\thalf, f\times\chi)\ll (qp^2)^{\frac14+\epsilon}(p^{-\frac18} +q^{\frac14}p^{-\frac14}).
\end{align}
For any fixed $\delta>0$, this gives subconvexity in the range $q \le  p^{1-\delta}$,  and is particularly strong when $q\le p^\epsilon$, for then it is of Burgess quality ($\kappa=\frac{3}{16}$).  In \cite{khan}, we established the subconvexity bound 
\begin{align}
\label{khan-sub} L(\thalf, f\times\chi)\ll (qp^2)^{\frac14+\epsilon}(q^{\frac14}p^{-\frac12}+q^{-\frac14})
\end{align}
  for $p^\delta \le  q \le  p^{2-\delta}$. This bound is particularly strong for $p^{1-\epsilon}\le q \le p^{1+\epsilon}$, where it is of Weyl quality ($\kappa=\frac{1}{6}$). In the special case that $\chi$ is real, Petrow and Young \cite{petyou3} proved the bound
 \[
L(\thalf, f\times\chi)\ll  (qp^2)^{\frac14+\epsilon}q^{\frac{1}{12}}p^{-\frac16},
\]
which gives subconvexity for $q\le p^{2-\delta}$ and is of Weyl quality for $q\le p^\epsilon$. 

In this paper, we prove a new subconvexity bound which has the advantage of covering at once essentially all ranges of $p$ and $q$ (except when they are very small) and being particularly strong for $p^{2-\epsilon}\le q\le p^{2+\epsilon}$, which is a new feature. We aimed particularly to address this latter range, for this is the point at which \eqref{khan-sub} fails to give any subconvexity bound.
\begin{theorem}\label{main}
Let $f$ be a holomorphic newform with weight $k\ge2$, prime level $q$ and trivial nebentypus. Let $\chi$ be a primitive Dirichlet character of prime modulus $p$. Suppose that $(q,p)=1$. We have 
\begin{align*}
L(\thalf, f\times \chi)\ll_{k,\epsilon} (qp^2)^{\frac14+\epsilon} (p^{-\frac14}+q^{-\frac18+\frac{\vartheta}{4}}),
\end{align*}
where $\vartheta=\frac{7}{64}$ denotes current best exponent towards the $GL(2)$ Ramanujan Conjecture, due to Kim and Sarnak \cite{kimsar}.
\end{theorem}
\noindent On the Ramanujan Conjecture, for any fixed $\delta>0$ our bound gives level aspect subconvexity as long as $p\ge q^{\delta}$ and $q\ge p^{\delta}$, and in the subrange $p^{2-\epsilon}\le q\le p^{2+\epsilon}$ we get the Burgess bound $L(\thalf, f\times \chi)\ll (qp^2)^{\frac{3}{16}+\epsilon}$. This latter bound is the best our result can do.

While the subconvexity bound \eqref{khan-sub} of \cite{khan} was deduced from an estimate for the second moment of $L(\thalf, f\times \chi)$, averaging over the newforms $f$, in this paper the main result Theorem \ref{main} is deduced from a fourth moment estimate. This is deeper than our previous paper and requires the spectral theory of automorphic forms. 
\begin{theorem} \label{4th}  Keep the notation in the statement of Theorem \ref{main}. For $k\ge 2$, we have
\begin{align} 
\label{4thm}  \sum_{f\in B_k^*(q)} |L(\thalf, f\times\chi)|^4  \ll_{k,\epsilon}  (qp)^\epsilon( qp+p^2q^{\half+\vartheta})
\end{align}
for any $\epsilon>0$. 
\end{theorem}
\noindent This fourth moment estimate immediately yields Theorem \ref{main}. For comparison, we note that the large sieve would give the bound $O(q^{1+\epsilon}p^{2+\epsilon})$ for the left hand side of \eqref{4thm} (see Lemma \ref{afe-use} and Theorem \ref{largesievethm}). Of course, the best expected bound is essentially $O(q^{1+\epsilon})$ on the Lindel\"of Hypothesis.

The study of the fourth moment of $GL(2)$ automorphic $L$-functions in the level aspect has a long history. In their influential paper \cite{dfi}, Duke, Friedlander, and Iwaniec were the first to prove a non-trivial upper bound for an amplified fourth moment, and from this the first subconvexity bound for automorphic $L$-functions in the level aspect. Kowalski, Michel, and VanderKam \cite{kmv} went further to prove an asymptotic for a mollified fourth moment, from which they obtained nonvanishing results in this context. Blomer, Harcos, Michel \cite{bloharmic2} proved estimates for an amplified fourth moment involving more general forms, allowing for non-trivial nebentypus. More recently, new techniques were introduced, which avoid the so-called shifted convolution problem that had been central to the aforementioned works. Kiral and Young \cite{kiryou} proved an estimate for the fifth moment that is sharp on the Ramanujan Conjecture. Blomer and Khan \cite{blokha2} established a `reciprocity' relation for the twisted fourth moment, and from this an estimate for an amplified fourth moment where the amplifier can be taken so long that it yields the fifth moment. This was generalized to the case of non-trivial nebentypus by the same authors in \cite{blokha}. 

In this paper, we study the fourth moment in a different direction. Instead of trying to insert a longer amplifier or mollifier, our goal is to reduce the family size in the fourth moment. Indeed in \eqref{4thm}, the $L$-functions are associated to forms of level $qp^2$, but we average over only those which arise as lifts of forms of level $q$. Note that in the range $p^{2-\epsilon}\le q\le p^{2+\epsilon}$, where we get our strongest subconvexity bound, the log(conductor) to log(family size) ratio is 8.  The fourth moment in small families has also been studied for other types of $L$-functions. See \cite{petyou2} for Dirichlet $L$-functions and \cite{jutila} for automorphic $L$-functions in the spectral aspect.

The general plan of attack is to introduce a root number into $|L(\thalf, f\times\chi)|^4$ using the functional equation, following \cite{blokha}. We then express the $L$-values as Dirichlet series using approximate functional equations, and apply the Petersson trace formula to obain a sum of Kloosterman sums. We dualize the sum by applying Voronoi summation, obtaining a new sum of Kloosterman sums. To this we apply Kuznetsov's formula in reverse to return to a sum of automorphic objects. The final step is to bound this using the spectral large sieve. The term $O(p^{2+\epsilon}q^{\half+\vartheta+\epsilon})$ seen in \eqref{4thm} arises from this process. The other term $O(q^{1+\epsilon}p^{1+\epsilon})$ arises from a diagonal term earlier on in the calculation (see Lemma \ref{bad-diag}). The sketch just described is part of a familiar sort of strategy in a classical approach to `spectral reciprocity'. We have already mentioned a few works in this regard. See also \cite{andkir, humkha, blolimil} for a further sample of the growing literature on spectral reciprocity type problems.  From \cite{blokha2}, for example, one sees that the fourth moment of $L(\half, f)$ over a family of level $q$ forms is related to a fourth moment of level $1$ forms, times the Hecke eigenvalue at $q$. In our problem, we consider the fourth moment of $L(\half, f\times \chi)$ over a family of level $q$ forms. Although we do not develop a precise reciprocity relation, we will roughly see that on the dual side, we get a fourth moment over level $p^4$ forms with nebentypus $\chi^2$, times the Hecke eigenvalue at $q$.

Throughout, we follow the $\epsilon$-convention: that is, $\epsilon$ will always be positive number which can be taken as small as we like, but may differ from one occurence to another. All implied constants may depend on $\epsilon$ and $k$.

{\bf Acknowledgement.} I am grateful to Matthew P. Young for helpful comments.


\section{Background}

\subsection{Automorphic forms and $L$-functions}\label{autl} 
Most of the facts in this subsection may be found in \cite[Chapter 14]{iwakow}. Let $S_k(q)$ denote the space of holomorphic cusp forms of prime level $q$, weight $k$, and trivial nebentypus. Let $S_k^*(q)\subset S_k(q) $ denote the space of newforms. Every $f\in S_k(q)$ has a Fourier series expansion
\begin{align*}
f(z)=\sum_{n=1}^\infty \rho_f(n) n^{\frac{k-1}{2}} e(nz)
\end{align*}
for $\Im(z)>0$. Let $B_k(q)$ denote an orthogonal basis of $S_k(q)$, which contains a basis $B_k^*(q)$ of $S_k^*(q)$, normalized so that $\rho_f(1)=1$ for every $f\in B_k^*(q)$. Thus for $f\in B_k^*(q)$, the $n$-th Fourier coefficient equals the $n$-th Hecke eigenvalue.

Let $\chi$ be primitive Dirichlet character of prime modulus $p$. For $f\in B_k(q)$, define
\begin{align*}
(f\times \chi)(z)=\sum_{n=1}^\infty \chi(n) \rho_f(n) n^{\frac{k-1}{2}} e(nz).
\end{align*}
For $f\in B_k^\star(q)$, this is a newform of level $qp^2$, weight $k$, and nebentypus $\chi^2$.  The associated $L$-function is entire and for $\Re(s)>1$ equals
\begin{align*}
L(s,f\times \chi)=\sum_{n=1}^\infty \frac{\chi(n) \rho_f(n)}{ n^s}.
\end{align*}
This satisfies the functional equation
\begin{align}
\label{feq} \Lambda(s, f\times \chi):= \Big(\frac{p q^\half  }{2\pi}\Big)^s\Gamma(s+\tfrac{k-1}{2}) L(s,f\times \chi)= \varepsilon(f\times\chi) \Lambda(1-s, f\times \overline{\chi}),
\end{align}
where
\[
\varepsilon(f\times\chi) =-  i^k \rho_f(q)q^{\half}\chi(q) \frac{  \tau(\chi)^2}{p}
\]
is of modulus 1. Thus the analytic conductor at $s=\frac12$ is $kp^2q$, or just $p^2q$ in the level aspect.

We now consider more general automorphic forms. Let $N\ge 1$ be an integer, $\psi$ a Dirichlet character of modulus $N$, and let cond$(\psi$) denote the modulus of the primitive character which induces $\psi$. Suppose that $\psi$ is even and cond$(\psi$) is squarefree and odd, since this is the case which we will need. Let $B_k(N,\psi)$ denote an orthogonal basis of holomorphic cusp forms of weight $k$ and nebentypus $\psi$ with respect to $\Gamma_0(N)$. Let $B(N,\psi)$ denote an orthogonal basis of Maass cusp forms with nebentypus $\psi$ with respect to $\Gamma_0(N)$, with $t_g$ denoting the spectral parameter of $g\in B(N,\psi)$. It is conjectured that $t_g$ is real, but the possibility of purley imaginary $t_g$ with $|t_g|\in(0,\half)$ has not been disproven (these are called exceptional eigenvalues). Whenever we wish to work with orthonormal bases, we will divide by the $L^2$-norm $\| g\|$ of the forms. Let $E_{\mathfrak{c}}(\cdot, \half+it)$ denote the Eisenstein series of nebentypus $\psi$ associated to a cusp $\mathfrak{c}$ which is singular for $\psi$ (see \cite{you} for definitions). For $\sigma=
(\begin{smallmatrix}
* & * \\
c & d 
\end{smallmatrix})
\in SL_2(\mathbb{R})$, let $j(\sigma, z):=cz+d$. We write the Fourier expansions of these objects around a singular cusp $\mathfrak{a}$ associated to a scaling matrix $\sigma_\mathfrak{a}$ as follows:
\begin{align*}
&g(\sigma_\frak{a} z)j(\sigma_\frak{a}, z)^{-k}=\sum_{n=1}^{\infty} \rho_{g_\frak{a}}(n) n^{\frac{k-1}{2}} e(nz) \ \ \text{for } g\in B_k(N,\psi),\\
&g(\sigma_\frak{a}z)= \sum_{n\neq 0} \rho_{g_\frak{a}}(n) n^{-\frac{1}{2}} W_{0,it_f} (4\pi |n| y)e(nx) \ \ \text{for } g\in B(N,\psi),\\
&E_\mathfrak{c}(\sigma_\frak{a}z,\thalf+it)=c_{1,\mathfrak{c}}(t)y^{\half+it}+c_{2,\mathfrak{c}}(t)y^{\half-it}+\sum_{n\neq 0} \rho_{\frak{c}\frak{a}}(n,t) n^{-\frac{1}{2}} W_{0,it_f} (4\pi |n| y)e(nx),
\end{align*}
where 
\[
W_{0,it_f} (4\pi |n| y)= \sqrt{n \pi y}K_{it_f}(2\pi |n| y)
\]
is a Whittaker function and $K_{it}$ is a Bessel function. For $\frak{a}=\infty$, we simply write 
\[
\rho_{g_\frak{a}}(n)=\rho_{g}(n) \ \ \text{ and } \ \ \rho_{\frak{c}\frak{a}}(n,t)=\rho_{\frak{c}}(n,t),
\]
and this matches the notation seen above.

Let $B_k^*(N,\psi)\subset B_k(N,\psi)$ and $B^*(N,\psi)\subset B(N,\psi)$ denote orthogonal bases of the space of newforms, normalized so that $\rho_g(1)=1$. For $g$ in $B_k^*(N,\psi)$ or $B^*(N,\psi)$, we have the Hecke multiplicativity relations
\[
\rho_g(n)\rho_g(m) =\sum_{\substack{d|(n,m)}} \psi(d) \rho_g\Big(\frac{nm}{d^2}\Big), \ \ \ \rho_g(nm)=\sum_{d|(n,m)} \mu(d) \psi(d) \rho_g\Big(\frac{n}{d}\Big)\rho_g\Big(\frac{m}{d}\Big).
\]
For $g$ in $B_k^*(N,\psi)$, we have the Ramanujan bound (due to Deligne),
\[
\rho_g(n)\ll n^\epsilon,
\]
and for $g\in B^*(N,\psi)$, we have the Kim-Sarnak bound \cite{kimsar}
\[
\rho_g(n)\ll n^{\vartheta+\epsilon},
\]
where $\vartheta=\frac{7}{64}$. For $g\in B_k^*(N)$ and $r$ dividing the square-free part of $N$, we have \cite[Theorem 2]{ogg}
\[
|\rho_g(r) | = r^{-\half}.
\]
Thus for $f\in B_k^*(q)$, we have $\rho_f(q)=\pm q^{-\half}$.

In order to appeal to the properties above, which are restricted to newforms, it will be useful to have bases $B_k(N,\psi)$ and $B(N,\psi)$ expressible in terms of lifts of newforms.

\begin{lemma} \label{specialbasis} We have an orthogonal basis
\begin{align}
B_k(N,\psi)= \bigcup_{\substack{\delta \ell | N }} \{ g^{(\delta)}: g\in B_k^*(\ell, \psi) \},
\end{align}
where
\[
g^{(\delta)}(z)=\sum_{r|\delta} \nu_\delta(r) r^\half g(rz)
\]
for some complex numbers $\nu_\delta(r)$ depending on $r, \delta$, and $g$, such that $ \nu_\delta(r)\ll (r\delta)^\epsilon$, where the understanding is that $B_k^*(\ell, \psi)$ is empty if $\psi$ is not a character mod $\ell$. Further, this basis is orthonormal if every $g\in B_k^*(\ell, \psi)$ for $\ell|N$ is $L^2$-normalized with respect to $\Gamma_0(N)$.

The same type of statement with obvious notational modification holds for $B(N,\psi)$.
\end{lemma}
\proof
See \cite[Proposition 7.1]{pet} and \cite[Section 6.6]{petyou2}.
\endproof
We do not reproduce it here, but analogous to Lemma \ref{specialbasis}, for the continuous spectrum we have an `orthonormal' basis involving newform Eisenstein series, whose Fourier coefficients satisfy the Hecke relations. This theory was developed in \cite{you}. See the remarks following \cite[Lemma 6.11]{petyou2}.


\subsection{Petersson trace formula} 
We state the Petersson trace formula.
\begin{theorem}
We have
\begin{align}
\label{trace} \sum_{f\in B_k(N)} \frac{\Gamma(k-1)}{(4\pi)^{k-1}\langle f, f\rangle } \rho_f(n_1)\overline{\rho_f(n_2)}=\delta_{(n_1=n_2)}+2\pi  i^{-k}\sum_{c\ge 1} \frac{S(n_1,n_2,cN)}{cN} J_{k-1}\Big(\frac{4\pi \sqrt{n_1n_2}}{cN}\Big),
\end{align}
where $S(n_1,n_2,cN)$ is the Kloosterman sum, $\langle f , f \rangle$ is the Petersson inner product for $\Gamma_0(N)$, $J_{k-1}(x)$ is a $J$-Bessel function, and $\delta_P$ equals $1$ if the statement $P$ is true, and $0$ otherwise. 
\end{theorem}
\proof
See \cite[Corollary 14.23]{iwakow}.
\endproof
By \cite[equation (2.24)]{iwamic}, for $N$ square-free and $f\in  B^*_k(N)$, we have
\begin{align}
\label{fnorm-ident} \frac{\Gamma(k-1)}{(4\pi)^{k-1}\langle f, f\rangle }=\frac{2\pi^2}{N(k-1)L(1,\sym^2 f)},
\end{align}
where $L(s,\sym^2 f)$ is the associated symmetric square $L$-function and $N^{-\epsilon}\ll  L(1,\sym^2 f) \ll  N^{\epsilon}$
by \cite{hofloc}.


\subsection{Spectral large sieve}\label{specsieve}

We state the spectral large sieve inequality.
\begin{theorem}\label{largesievethm}
For any sequence $\{\alpha_n\}$ of complex numbers, we have 
\begin{align*}
&\sum_{\substack{k\le T\\ k\equiv 0\bmod 2}} \frac{(k-1)!}{(4\pi)^{k-1}}  \sum_{g\in B_k(N,\psi)}\frac{1}{\| g \|^2} \Big| \sum_{M\le m \le 2M} \alpha_m \rho_{g_\mathfrak{a}}(m )\Big|^2 \ll \Big(T^2+ \frac{M^{1+\epsilon}}{N}\Big)\| \alpha \|^2,\\
& \sum_{\substack{g\in B(N,\psi)\\ t_f\le T}} \frac{1}{\| g \|^2}  \frac{1}{\cosh(\pi t_g)} \Big| \sum_{M\le m \le 2M} \alpha_m \rho_{g_\frak{a}}(m)\Big|^2 \ll \Big(T^2+ \frac{M^{1+\epsilon}}{N}\Big)\| \alpha \|^2,\\
& \sum_{\mathfrak{c} \text{ singular} } \intt \frac{1}{\cosh(\pi t)} \Big| \sum_{M\le m \le 2M} \alpha_m \rho_{\mathfrak{c}\frak{a}}(m,t)\Big|^2 \ll \Big(T^2+ \frac{M^{1+\epsilon}}{N}\Big)\| \alpha \|^2,
\end{align*}
where $\| g\|^2$ is the Petersson inner product of $g$ with itself and $\| \alpha \|^2=\sum_{M\le m \le 2M} |\alpha_m|^2 $.
\end{theorem}
\proof
This is almost given by \cite[Proposition 4.7]{drap}, except that there, the bounds on the right hand side equal $(T^2+ \mathrm{cond}(\psi)^\half \frac{M^{1+\epsilon}}{N})\| a \|^2$. Under the assumption that $\mathrm{cond}(\psi)$ is squarefree and odd, it was observed in \cite[Proposition 2.5]{zac} that this extra factor $\mathrm{cond}(\psi)$ can be removed.
\endproof

When the sums on the left hand side of the large sieve inequality are restricted to newforms, normalized to have first Fourier coefficient equal to $1$, we can use
\begin{align}
\label{obs} N^{1-\epsilon}k^{-\epsilon} \frac{(k-1)!}{(4\pi)^{k-1}} \ll \| g \|^2\ll N^{1+\epsilon}k^\epsilon \frac{(k-1)!}{(4\pi)^{k-1}}
\end{align}
for $g$ holomorphic, and 
\[
N^{1-\epsilon}(1+|t_g|)^{-\epsilon} \frac{1}{\cosh(\pi t_g)} \ll \| g \|^2\ll N^{1+\epsilon}(1+|t_g|)^{\epsilon}\frac{1}{\cosh(\pi t_g)}
\]
for $g$ a Maass form. See \cite[Equation (2.3)]{petyou1}. For analogous bounds for newform Eisenstein series, see \cite[Equation (2.7)]{petyou1}.


\subsection{Kuznetsov's trace formula}
Follow the notation from the previous section.
Let $\phi$ be a smooth function compactly supported on the positive reals. Define the following transforms of $\phi$ against Bessel functions:
\begin{align*}
&\dot{\phi}(k)=\frac{i^k}{\pi}\int_0^\infty J_{k-1}(x)\phi(x)\frac{dx}{x},\\
&\tilde{\phi}(t)=\frac{2\pi i}{\sinh(\pi t)}\int_0^\infty (J_{2it}(x)-J_{-2it}(x)) \phi(x) \frac{dx}{x},\\
&\check{\phi}(t)= 8\cosh(\pi t)\int_0^\infty K_{2it}(x) \phi(x) \frac{dx}{x}.\\
\end{align*}
We have the following estimates for these transforms.
\begin{lemma} \label{transf-bounds}
(a) If $\phi(x)$ is supported on $0<X<x<2X$ and satisfies $\phi^{(j)}(x)\ll_j X^{-j}$ for all $j\ge 0$, then for $t\in\mathbb{R}$ we have
\begin{align*}
\dot{\phi}(t), \tilde{\phi}(t), \check{\phi}(t)\ll_l \frac{1+|\log X|}{1+X}\Big(\frac{1+X}{1+|t|}\Big)^l
\end{align*}
for any $l\ge 0$. 

(b) If $\phi(x)$ is supported on $0<X<x<2X$ and satisfies $\phi^{(j)}(x)\ll_j (\frac{X}{Z})^{-j}$ for all $j\ge 0$, then for $t\in(-\frac{i}{4}, \frac{i}{4})$, we have
\begin{align*}
\tilde{\phi}(t), \check{\phi}(t)\ll \frac{1+(\frac{X}{Z})^{-2|\Im t|}}{1+(\frac{X}{Z})}.
\end{align*}

(c) If $\phi(x)=e^{ax}\psi(x)$, where $a=\pm 1$, and $\psi(x)$ is supported on $1\le X<x<2X$ and satisfies $\phi^{(j)}(x)\ll_j X^{-j}$ for all $j\ge 0$, then for $t\in\mathbb{R}$, we have
\begin{align*}
\dot{\phi}(t),  \tilde{\phi}(t) \ll X^{-\half+\epsilon}.
\end{align*}
For $|t|> X^{\half+\epsilon}$, we have 
\[
\dot{\phi}(t), \tilde{\phi}(t) \ll_l (|t|+X)^{-l}X^\epsilon
\]
for any $l\ge 0$.
\end{lemma}
\proof
For parts (a) and (b), see \cite[Lemma 1]{bloharmic2}, and for part (c), see \cite[Lemma 3, Remark 1]{jut2}.
\endproof

Given two singular cusps $\mathfrak{a}$ and $\mathfrak{b}$ of $\Gamma_0(N)$, define 
\begin{align*}
&\mathcal{H}=\sum_{\substack{k>0\\ k\equiv 0 \bmod 2}} \sum_{g\in B_k(N,\psi)}  \dot{\phi}(k) \frac{(k-1)!}{(4\pi)^{k-1}} \frac{\rho_{g_\frak{a}}(n)\overline{\rho_{g_\frak{a}}(m)}}{\| g\|^2}\\
&\mathcal{M}= \sum_{g\in B(N,\psi)} \frac{1}{\cosh(t_g)} \tilde{\phi}(t_g) \frac{\rho_{g_\frak{a}}(n)\overline{\rho_{g_\frak{a}}(m)}}{\| g\|^2}, \ \ \ \\
&\mathcal{E}=\sum_{\frak{c}} \intt  \frac{1}{4\pi \cosh(t)} \tilde{\phi}(t) \rho_{\frak{a}\frak{c}}(n)\overline{\rho_{\frak{c}\frak{a}}(m)}\\
&\mathcal{M}^\prime=\sum_{g\in B(N,\psi)} \frac{1}{\cosh(t_g)} \check{\phi}(t_g) \frac{\rho_{g_\frak{a}}(-n)\overline{\rho_{g_\frak{a}}(m)}}{\| g\|^2}\\
&\mathcal{E}^\prime=\sum_{\frak{c}} \intt  \frac{1}{4\pi \cosh(t)} \check{\phi}(t) \rho_{\frak{a}\frak{c}}(-n)\overline{\rho_{\frak{c}\frak{a}}(m)}.
\end{align*}
We now state a special case of Kuznetsov's theorem which we will need.
\begin{theorem} \label{kuznetsov} Keep the notation of section \ref{specsieve}. Suppose that the level $N=rs$, where the modulus of $\psi$ divides $r$, and $(r,s)=1$. Let $\mathfrak{a}=\infty$ and $\frak{b}=1/s$ (these are singular cusps). Let $\phi$ be a smooth function compactly supported on the positive reals. For $n,m\ge 1$ we have
\begin{align*}
&\sum_{\substack{c\ge 1\\ (c,r)=1}} \frac{1}{cs\sqrt{r}} \overline{\psi}(c) e\Big(\frac{m\overline{s}}{r}\Big) S(m\overline{r},n, cs)  \phi\Big(\frac{4\pi\sqrt{nm}}{cs\sqrt{r}}\Big)=\mathcal{H}+\mathcal{M}+\mathcal{E},\\
&\sum_{\substack{c\ge 1\\ (c,r)=1}} \frac{1}{cs\sqrt{r}} \overline{\psi}(c) e\Big(\frac{m\overline{s}}{r}\Big) S(-m\overline{r},n, cs) \phi\Big(\frac{4\pi\sqrt{nm}}{cs\sqrt{r}}\Big) =\mathcal{M}^\prime+\mathcal{E}^\prime,\\
\end{align*}
where the Fourier coefficients in $\mathcal{H},\mathcal{M},\mathcal{E},\mathcal{M}^\prime,\mathcal{E}^\prime$ are determined by a suitable choice of scaling matrices.
\end{theorem}
\proof
See \cite[Lemmas 4.3 and 4.5]{drap}. The first of these lemmas indicates the scaling matrices that we need.
\endproof


\subsection{Bessel functions}

In this section we collect various estimates for Bessel functions.
\begin{lemma} \label{jlemma}
For $(qp)^{-100}\le x\le (pq)^{\epsilon}$, we have 
\[
J_{k-1}(4\pi x)=xW_1(x)
\] 
where $W_1(x)$ is a smooth function satisfying $W_1^{(j)}(x)\ll_j (qp)^\epsilon$ for all $j\ge 0$.
For $x> (qp)^{\epsilon}$, we have
\[
J_{k-1}(4\pi x)=\Re\Big( \frac{e(2x)}{\sqrt{x}}W_2(x) \Big),
\]
where $W_2$ is a smooth function satisfying
\begin{align*}
 x^jW_2^{(j)}(x)\ll_{j} 1
\end{align*}
for all $j\ge 0$. 
\end{lemma}

\proof 

When $x\le (qp)^{-\epsilon}$, use the power series \cite[equation 8.440]{gr},
\[
J_{k-1}(x)=\sum_{\ell=0}^\infty \frac{(-1)^\ell (x/2)^{k-1+2\ell}}{\ell! (k+\ell)!}.
\]
When $x>(qp)^{-\epsilon}$, use the integral representation \cite[page 206]{wat}
\[
J_{k-1}(4\pi x)=  \Re\Bigg(\frac{e(2x)}{\sqrt{x}} \frac{e(-\frac{k}{4}+\frac{1}{8})}{\pi\sqrt{2}\Gamma(k-\half)}\int_0^\infty e^{-u}u^{k-\frac32}\Big(1+\frac{iu}{8\pi x}\Big)^{k-\frac32} du\Bigg).
\]

\endproof 

\begin{lemma} \label{jlem}
(a) For $x\ge (qp)^{\epsilon}$, we have $K_0(x)\ll (qp)^{-100}$, and for $(qp)^{-100} < x< (qp)^{\epsilon}$, we have 
\[
K_0^{(j)} (x)\ll_j (qp)^\epsilon (1+x^{-j})
\]
for $j\ge 0$. 

(b) For $x\ge (qp)^{\epsilon}$, we have 
\begin{align}
\label{y0long} Y_0(x)=\Im \Big(\frac{e^{ix}}{\sqrt{x}} W(x)\Big),
\end{align}
for some smooth function $W$ with $x^{-j} W^{(j)}(x)\ll_j 1$ for $j\ge 0$.
For $(qp)^{-100} < x< (qp)^{\epsilon}$, we have 
\[
Y_0^{(j)} (x)\ll_j (qp)^\epsilon(1+x^{-j})
\]
for $j\ge 0$.
\end{lemma}
\proof
(a) By \cite[page 208]{wat}, we have
\[
K_0(x)=e^{-x} 2^{-\half} \int_0^\infty e^{-ux} u^{-\half}(1+\thalf u)^{-\half} du.
\]
Suppose first that $x\ge (qp)^{\epsilon}$. In this range, the integral is $O(1)$ and due to the exponential decay of the $e^{-x}$ factor in front of the integral, we we get the desired bound $K_0(x)\ll (qp)^{-100}$.  Now suppose $(qp)^{-100} < x< (qp)^{\epsilon}$. To see the claim in this range, observe that
\[
 \int_0^\infty e^{-ux} u^{-\half}(1+\thalf u)^{-\half} du \ll  \int_0^1 u^{-\half} du + \int_1^{\max(1,x^{-2})} u^{-1}du  +\int_{\max(1,x^{-2})}^\infty e^{-ux} du\ll (qp)^\epsilon,
\]
and 
\[
\frac{d^j}{dx^j}  \int_0^\infty e^{-ux} u^{-\half}(1+\thalf u)^{-\half} du \ll  \int_0^\infty u^{j-1} e^{-ux} du \ll x^{-j}.
\]
for $j\ge 1$.

(b) By \cite[page 206]{wat}, we have
\begin{align}
\label{ybes-big} Y_0(x)=\Im \Big(\frac{e^{ix}}{\sqrt{x}} W(x)\Big),
\end{align}
where
\[
W(x)= \frac{\sqrt{2}}{\pi} \int_0^\infty e^{-u} u^{-\half}\Big(1+\frac{iu}{2x}\Big)^{-\half} du.
\]
For $x\ge (qp)^\epsilon$, we have $W^{(j)}(x)\ll  x^{-j}$ for $j\ge 0$. Now consider the range $(qp)^{-100} < x< (qp)^{\epsilon}$. For the subrange $(qp)^{-\epsilon} < x< (qp)^{\epsilon}$, we use the integral representation above to see that $\frac{d^j}{dx^j} Y_0(x)\ll_j (qp)^\epsilon$. For $(qp)^{-100} < x< (qp)^{-\epsilon}$, we use the power series \cite[equations 8.402, 8.403.2]{gr},
\[
Y_0(x)=\frac{2}{\pi}\Big(\ln(\thalf x)+\gamma\Big)\sum_{l=0}^\infty \frac{(-x^2)^l}{4^l(l!)^2}+\frac{2}{\pi}\sum_{l=0}^\infty  \frac{(-x^2)^l}{4^l(l!)^2} \sum_{i=1}^{l+1} \frac{1}{i},
\]
to get that $Y_0^{(j)}(x)\ll_j (qp)^\epsilon(1+x^{-j})$ for $j\ge 0$.
\endproof


\subsection{Voronoi summation}

In this section we state some summation formulae for divisor type functions. Define
\[
\tau_{\overline{\chi}^2}(n)=\sum_{l|n}\overline{\chi}^2 (l), \ \ \ \ \ \tau(n)=\sum_{l|n} 1.
\]

\begin{lemma} \label{vor} Let $g$ be a smooth function compactly supported on the positive reals. Let $\chi$ be a primitive Dirichlet character of modulus $p$, and $c$, $b$ integers with $(c,pb)=1$. We have
\begin{align}
\label{vor1} \sum_{n\ge 1} \chi(n) \tau(n)e\Big(\frac{nb}{c}\Big)g(n)&=-\frac{2\pi }{cp} \frac{\tau(\chi)}{\tau(\overline{\chi} )} \sum_{n\ge 1} \chi(c^2) \overline{\chi}(-n) \tau(n) e\Big(\frac{-n \overline{p}^2\overline{b}}{c}\Big)\int_0^\infty Y_0\Big(\frac{4\pi\sqrt{nx}}{cp}\Big) g(x) dx\\
\nonumber &+\frac{4}{cp}\frac{\tau(\chi)}{\tau(\overline{\chi} )} \sum_{n\ge 1}  \chi(c^2) \overline{\chi}(-n)  \tau(n) e\Big(\frac{n \overline{p}^2\overline{b}}{c}\Big) \int_0^\infty K_0\Big(\frac{4\pi\sqrt{nx}}{cp}\Big) g(x) dx,\\
\label{vor2} \sum_{n\ge 1} \chi(n)\tau_{\overline{\chi}^2}(n) e\Big(\frac{nb}{c}\Big)g(n)=& -\frac{2\pi }{cp} \sum_{n\ge 1}\chi(-n) \tau_{\overline{\chi}^2}(n) e\Big(\frac{-n \overline{p}^2\overline{b}}{c}\Big)\int_0^\infty Y_0\Big(\frac{4\pi\sqrt{nx}}{cp}\Big) g(x) dx\\
\nonumber &+\frac{4}{cp} \sum_{n\ge 1} \chi(n) \tau_{\overline{\chi}^2}(n) e\Big(\frac{n \overline{p}^2\overline{b}}{c}\Big) \int_0^\infty K_0\Big(\frac{4\pi\sqrt{nx}}{cp}\Big) g(x) dx,
\end{align}
where $Y_0$ and $K_0$ denote Bessel functions. 
\end{lemma}
\proof
We first consider \eqref{vor1}. Using \cite[equation (3.12)]{iwakow}, we have
\begin{align}
\label{vor1-1} \sum_{n\ge 1} \chi(n) \tau(n)e\Big(\frac{nb}{c}\Big)g(n)&= \frac{1}{\tau(\overline{\chi})} \sum_{a\bmod p} \overline{\chi}(a) \sum_{n\ge 1}  \tau(n)e\Big(\frac{na}{p}+\frac{nb}{c}\Big)g(n)\\
\nonumber &=\frac{1}{\tau(\overline{\chi})} \sum_{a\bmod p} \overline{\chi}(a) \sum_{n\ge 1}  \tau(n)e\Big(\frac{n(ac+bp)}{cp}\Big)g(n).
\end{align}
We can restrict to $(a,p)=1$ since otherwise $\overline{\chi}(a)=0$. Observe that the inverse of $ac+bp$ in $\mathbb{Z}_{cp}^*$ is $c(\overline{c})_p^2(\overline{a})_p+p(\overline{p})_c^2(\overline{b})_c$, where $(\overline{c})_p$ denotes the inverse of $c$ in $\mathbb{Z}_{p}^*$ , etc.  Now, applying the Voronoi summation formula for the divisor function \cite[equation (4.49)]{iwakow}, we get that \eqref{vor1-1} equals
\begin{align}
\label{vor1-2} &\frac{1}{\tau(\overline{\chi})} \summ_{a\bmod p} \overline{\chi}(a) \frac{1}{cp}\int_0^\infty (\log x + 2\gamma -2 \log(cp))g(x)dx\\
\nonumber &-\frac{2\pi }{cp} \frac{1}{\tau(\overline{\chi} )} \summ_{a\bmod p} \overline{\chi}(a) \sum_{n\ge 1} \tau(n) e\Big(\frac{-n(c(\overline{c})_p^2(\overline{a})_p+p(\overline{p})_c^2(\overline{b})_c)}{cp}\Big)\int_0^\infty Y_0\Big(\frac{4\pi\sqrt{nx}}{cp}\Big) g(x) dx\\
\nonumber &+\frac{4}{cp}\frac{1}{\tau(\overline{\chi} )} \summ_{a\bmod p} \overline{\chi}(a) \ \sum_{n\ge 1} \tau(n) e\Big(\frac{n(c(\overline{c})_p^2(\overline{a})_p+p(\overline{p})_c^2(\overline{b})_c)}{cp}\Big)\int_0^\infty K_0\Big(\frac{4\pi\sqrt{nx}}{cp}\Big) g(x) dx.
\end{align}
We write
\[
e\Big(\frac{-n(c(\overline{c})_p^2(\overline{a})_p+p(\overline{p})_c^2(\overline{b})_c)}{cp}\Big)=e\Big(\frac{-n \overline{c}^2 \overline{a}}{p}\Big)e\Big(\frac{-n \overline{p}^2 \overline{b}}{c}\Big)
\]
and evaluate the $a$-sums in \eqref{vor1-2}. By \cite[equation (3.9)]{iwakow}, we get that the first sum vanishes, and by using \cite[equation (3.12)]{iwakow} for the other sums, we get the required formula \eqref{vor1}.

Now we turn to \eqref{vor2}. If $\chi$ is quadratic, then $\tau_{\overline{\chi}^2}(n)=\tau(n)$ for $(n,p)=1$, so we can refer to \eqref{vor1}. If $\chi$ is not quadratic, then \eqref{vor2} is justified the same way as \eqref{vor1}, except that we use the summation formula \cite[equation (4.67)]{iwakow} for $\tau_{\overline{\chi}^2}(n)$ instead of the Voronoi summation formula for $\tau(n)$.
\endproof


\subsection{Approximate functional equation}\label{afe-section}

In this section, we get a handle on the fourth power of the central $L$-value by expressing it in terms of a Dirichlet series. The main point to note is the introduction of the root number (from the functional equation) into this expression, following \cite{blokha}. This way, when we later apply Voronoi summation to the $n$ and $m$ sums below, it will not be circular.

\begin{lemma} \label{afe-use}
Let $f\in B_k^*(q)$. We have that $|L(\half, f\times \chi)|^4$ equals a linear combination of 
\begin{align}
\label{expression}q^\frac12  \sum_{\substack{d_1,d_2\ge 1\\ (d_1d_2,pq)=1}} \frac{\chi^2(d_2)}{d_1d_2}  \sum_{\substack{n,m \ge 1\\ (m,q)=1}}  \frac{ \chi(n) \tau_{\overline{\chi}^2}(n)  \rho_f(nq) }{n^\half}\frac{  \chi(m) \tau(m)  \overline{\rho_f(m)}
 }{m^\half}V\Big(\frac{nd_1^2}{qp^2}\Big)V\Big(\frac{md_2^2}{qp^2}\Big) +O(q^{-\half+\epsilon}p^{2+\epsilon}),
  \end{align}
and the same sum with $\chi(m)$ replaced by $\overline{\chi}(m)$ and $\chi^2(d_2)$ replaced by $\overline{\chi}^2(d_2)$, where
\[
V(x)=\frac{1}{2\pi i} \int_{(2)} (4\pi^2 x)^{-s}  \frac{\Gamma(\frac{k}{2}+s)^2}{\Gamma(\frac{k}{2})^2} \frac{ds}{s}
\] 
for $x>0$, and the constants involved in the linear combination are of modulus 1 and independent of $f$. 
\end{lemma}
\proof 
For $f\in B_k^*(q)$, we have by the functional equation \eqref{feq}, that
\[
 |L(\thalf, f\times \chi)|^4=|L(\thalf, f\times \chi)|^2L(\thalf, f\times \chi)L(\thalf, f\times \overline{\chi})= \varepsilon(f\times\overline{\chi}) |L(\thalf, f\times \chi)|^2 L(\thalf,f\times \chi)^2.
\]
Using the functional equations of $|L(s, f\times \chi)|^2$ and $L(s, f\times \chi)^2$, we can express their central values by approximate functional equations \cite[Theorem 5.3]{iwakow} as follows:
\[
|L(\thalf, f\times \chi)|^2= 2\sum_{n_1,n_2\ge 1} \frac{\rho_f(n_1)\chi(n_1)\rho_f(n_2)\overline{\chi}(n_2)}{(n_1n_2)^\half} V\Big(\frac{n_1n_2}{qp^2}\Big)
\]
and
\begin{align*}
L(\thalf, f\times \chi)^2= &\sum_{m_1,m_2\ge 1} \frac{\rho_f(m_1)\chi(m_1)\rho_f(m_2)\chi(m_2)}{(m_1m_2)^\half} V\Big(\frac{m_1m_2}{qp^2}\Big) \\
+ &\varepsilon(f\times \chi)^2\sum_{m_1,m_2\ge 1} \frac{\rho_f(m_1)\overline{\chi}(m_1)\rho_f(m_2)\overline{\chi}(m_2)}{(m_1m_2)^\half} V\Big(\frac{m_1m_2}{qp^2}\Big),
\end{align*}
where $V(x)$ is a weight function as defined in the lemma, which effectively restricts the sums to $n_1n_2\le q^{1+\epsilon}p^{2+\epsilon}$ and  $m_1m_2\le q^{1+\epsilon}p^{2+\epsilon}$. For $x\gg (qp)^{-100}$, have $x^j V^{(j)}(x) \ll_j (qp)^\epsilon$ for all $j\ge 0$ by moving the line of integration to $\Re(u)=\epsilon$. Note that since $\rho_f(q)q^{-\half}=\pm 1$,
we have that $\varepsilon(f\times \chi)^2$ is independent of $f$. Thus the fourth power $|L(\half, f\times \chi)|^4$ equals a constant of modulus 1, depending on $\chi, p$, and $q$, times 
\[
\rho_f(q)q^{\half}  \sum_{\substack{n_1,n_2\ge 1\\m_1,m_2\ge 1}} \frac{\rho_f(n_1)\chi(n_1)\rho_f(n_2)\overline{\chi}(n_2)}{(n_1n_2)^\half} V\Big(\frac{n_1n_2}{qp^2}\Big)  \frac{\rho_f(m_1)\chi(m_1)\rho_f(m_2)\chi(m_2)}{(m_1m_2)^\half} V\Big(\frac{m_1m_2}{qp^2}\Big),
\]
plus a similar sum with $\chi(m_1),\chi(m_2)$ replaced with $\overline{\chi}(m_1),\overline{\chi}(m_2)$. 
Using the Hecke multiplicativity relations for $f\in B_k^*(q)$, we have
\[
 \rho_f(q) \rho_f(n_1)\rho_f(n_2)=\sum_{\substack{d_1|(n_1,n_2)\\(d_1,q)=1}} \rho_f\Big(\frac{qn_1n_2}{d_1^2}\Big).
\]
We replace $n_1,n_2$ by $n_1d_1, n_2d_1$ and then write $n=n_1n_2$, and similarly for $m_1, m_2$. We observe that
\[
\sum_{n_1n_2=n}\chi(n_1)\overline{\chi}(n_2)=\chi(n)\sum_{n_2|n} \overline{\chi}^2(n_2),
\]
and similarly for the $m_1, m_2$ sum. Writing $\rho_f(m)=\overline{\rho_f(m)}$ since the Fourier coefficients are anyway real, we obtain the main term of the lemma without the condition $(m,q)=1$. If we impose this extra condition then we must separately consider the contribution of the terms with $m=m'q$. Using that $\rho_g(nq)$ and $\overline{\rho_g(m'q)}$ are each bounded by $q^{-\half+\epsilon}$, we obtain by trivially summing over $n$ and $m'$ that this separate contribution is bounded by the error displayed in the lemma.
\endproof


 \subsection{Mellin inversion}
 We will use Mellin inversion to separate variables in the argument of a function.
 \begin{lemma}\label{mellin-split}
Let $W(x)$ be a real valued smooth function compactly supported on $(qp)^{-100}<x<\infty$, satisfying $x^j W^{(j)}(x)\ll_j (qp)^\epsilon $ for $j\ge 0$. Let $\tilde{W}(s)$ denote the Mellin transform.
We have
\[
W(x)=\frac{1}{2\pi i}\int_{-i(qp)^\epsilon}^{i(qp)^\epsilon} x^{-s} \tilde{W}(s) ds +O((qp)^{-100}),
\]
where $\tilde{W}(s)\ll_j (qp)^\epsilon$ within the limits of the integral.
\end{lemma}
\proof
The Mellin transform
\[
\tilde{W}(s)=\int_0^\infty W(x) x^{s-1} dx
\]
exists for all $s\in\mathbb{C}$ and satisfies $\tilde{W}(it)\ll_j (qp)^\epsilon (1+|t|)^{-j}$ for all $j\ge 0$, by bounding trivially for $|t|\le 1$ and integrating by parts with respect to $x$ for $|t|> 1$. By Mellin inversion, we have
\[
W(x)=\frac{1}{2\pi i}\int_{-i\infty}^{i\infty} x^{-s} \tilde{W}(s) ds.
\]
The bound given above for $\tilde{W}(it)$ can be used to restrict the integral to $|t|\le (qp)^\epsilon$ up to an error of $O((qp)^{-100})$.
\endproof


\section{Preliminary work}

The goal in this section is to reduce Theorem \ref{4th} to a problem on estimating a sum of Kloosterman sums, which will be our main object of study.

\begin{lemma} \label{suffices}
To prove Theorem \ref{4th}, it suffices to establish
\begin{multline}
\label{need02}  \frac{1}{q^{\frac12}}  \sum_{\substack{n,m,c \ge 1\\  (c,pq)=1}} \frac{ \chi(n) \tau_{\overline{\chi}^2}(n)  }{n^\half}\frac{  \chi(m) \tau(m)   }{m^\half}\frac{S(n,m\overline{q},c)}{c}  J_{k-1}\Big(\frac{4\pi \sqrt{nm}}{c\sqrt{q}}\Big) V\Big(\frac{nd_1^2}{qp^2}\Big)  V\Big(\frac{md_2^2}{qp^2}\Big)\\
 \ll (qp)^\epsilon p^2q^{-\half+\vartheta},
\end{multline}
for any $d_1,d_2\le (qp)^{100}$, and the same bound with $\chi(m)$ replaced by $\overline{\chi}(m)$.
\end{lemma}
\proof
Using \eqref{fnorm-ident}, for Theorem \ref{4th} we need to prove 
\begin{align*}
\sum_{f\in B^*_k(q)}  \frac{\Gamma(k-1)}{(4\pi)^{k-1}\langle f, f\rangle }  | L(\thalf, f\times \chi)|^4  \ll  (qp)^\epsilon( p+p^2q^{-\half+\vartheta}).
\end{align*}
We write $|L(\thalf, f\times \chi)|^4$ using Lemma \ref{afe-use}. Then since  $|L(\thalf, f\times \chi)|^4$ is positive, we may enlarge $B^*_k(q)$ to $B_k(q)$. Thus it suffices to prove
\begin{multline*}
q^\frac12  \sum_{\substack{d_1,d_2\ge 1\\ (d_1d_2,pq)=1}} \frac{\chi^2(d_2)}{d_1d_2}  \sum_{\substack{n,m \ge 1\\ (m,q)=1}}  \frac{ \chi(n) \tau_{\overline{\chi}^2}(n)  }{n^\half}\frac{  \chi(m) \tau(m)   }{m^\half}V\Big(\frac{nd_1^2}{qp^2}\Big)V\Big(\frac{md_2^2}{qp^2}\Big)\\ \times \sum_{f\in B_k(q)} \frac{\Gamma(k-1)}{(4\pi)^{k-1}\langle f, f\rangle } \rho_f(nq)\overline{\rho_f(m)}
 \ll (qp)^\epsilon( p+p^2q^{-\half+\vartheta}),
\end{multline*}
and the same bound with $\chi(m)$ replaced by $\overline{\chi}(m)$ and $\chi^2(d_2)$ replaced by $\overline{\chi}^2(d_2)$. Note that the error term above absorbs the contribution of the error term seen in Lemma \ref{afe-use}. We apply the Petersson trace formula to the innermost sum on the left hand side. The contribution of the $\delta_{(nq=m)}$ `diagonal' terms is zero because of the condition $(m,q)=1$. Thus we are left to consider the `off-diagonal' terms, for which it suffices to prove, after forsaking cancellation in the $d_1,d_2$-sums, that 
\begin{multline}
\label{need1}  q^\frac12  \sum_{\substack{n,m,c \ge 1\\ (m,q)=1}} \frac{ \chi(n) \tau_{\overline{\chi}^2}(n)  }{n^\half}\frac{  \chi(m) \tau(m)   }{m^\half}\frac{S(nq,m,cq)}{cq}  J_{k-1}\Big(\frac{4\pi \sqrt{nmq}}{cq}\Big) V\Big(\frac{nd_1^2}{qp^2}\Big)  V\Big(\frac{md_2^2}{qp^2}\Big)\\
 \ll (qp)^\epsilon( p+p^2q^{-\half+\vartheta})
\end{multline}
for any $d_1,d_2\le (qp)^{100}$. Next we wish to restrict to $(c,p)=1$, but we relegate most of the work to a separate lemma. The contribution of the terms in \eqref{need1} with $c=c'p$, using the Petersson trace formula to sum over $c'\ge 1$, equals 
\begin{align*}
q^\frac12  \sum_{\substack{n,m \ge 1\\ (m,q)=1}}  \frac{ \chi(n) \tau_{\overline{\chi}^2}(n)  }{n^\half}\frac{  \chi(m) \tau(m)   }{m^\half}V\Big(\frac{nd_1^2}{qp^2}\Big)V\Big(\frac{md_2^2}{qp^2}\Big) \sum_{g\in B_k(pq)} \frac{\Gamma(k-1)}{(4\pi)^{k-1}\langle g, g\rangle } \rho_g(nq)\overline{\rho_g(m)}.
 \end{align*}
This satisfies the bound $O(q^\epsilon p^{1+\epsilon})$ using Lemma \ref{bad-diag} below, which is an admissible error in \eqref{need1}. We now simplify the Kloosterman sum in \eqref{need1}. Write $c=c'q^ l $ with $(c',q)=1$. If $ l \ge1$, we claim that $S(nq,m,cq)=0$. To see this, oberve by \cite[equation (1.59)]{iwakow} that
\[
 S(nq,m,cq)=S(nq,m,c'q^{ l +1})=S(n\overline{q}^{ l }, m, \overline{q}^{ l +1},c')S(nq\overline{c'},m\overline{c'},q^{ l +1})
\]
and
\begin{align*}
\label{kloo} S(nq\overline{c'},m\overline{c'},q^{ l +1})&=\sum_{\substack{0 \le a< q^ l \\ (a,q)=1}} \sum_{\substack{0 \le b < q }}  e\Big(\frac{nq\overline{c'}(a+bq^ l ) + m\overline{c'}(\overline{a+bq^ l })}{q^{ l +1}}\Big)\\
&= \sum_{\substack{0 \le a< q^ l \\ (a,q)=1}} \sum_{\substack{0 \le b < q }}   e\Big(\frac{nq\overline{c'}a+ m\overline{c'}(\overline{a}-\overline{a}^2 bq^ l )}{q^{ l +1}}\Big)=0,
\end{align*}
because the $b$-sum vanishes. This uses $(m,q)=1$. So we can assume that $(c,q)=1$, in which case we have
\[
 S(nq,m,cq)=S(n,m\overline{q},c)S(0,m\overline{c},q)=-S(n,m\overline{q},c).
\]
Thus to prove \eqref{need1}, it suffices to establish
\begin{multline}
\label{nee} q^{-\half}   \sum_{\substack{n,m,c \ge 1\\ (m,q)=1\\ (c,pq)=1}} \frac{ \chi(n) \tau_{\overline{\chi}^2}(n)  }{n^\half}\frac{  \chi(m) \tau(m)   }{m^\half}\frac{S(n,m\overline{q},c)}{c}  J_{k-1}\Big(\frac{4\pi \sqrt{nm}}{c\sqrt{q}}\Big) V\Big(\frac{nd_1^2}{qp^2}\Big)  V\Big(\frac{md_2^2}{qp^2}\Big)\\
 \ll (qp)^\epsilon p^2q^{-\half+\vartheta}.
\end{multline}
At this point we can remove the restriction $(m,q)=1$ because the terms with $m=m'q$ contribute, after using inclusion-exclusion to express the condition $(c,pq)=1$,
\begin{multline}
\label{remove}  q^{-1} \sum_{l|pq} \mu(l) \sum_{\substack{d_1,d_2\ge 1\\ (d_1d_2,pq)=1}} \frac{\chi^2(d_2)}{d_1d_2}  \sum_{\substack{n,m',c' \ge 1}} \frac{ \chi(n) \tau_{\overline{\chi}^2}(n)  }{n^\half}\frac{  \chi(m'q) \tau(m'q)   }{(m')^\half}\frac{S(n,m',c'l)}{c'l} \\ J_{k-1}\Big(\frac{4\pi \sqrt{nm'}}{c'l}\Big) V\Big(\frac{nd_1^2}{qp^2}\Big)  V\Big(\frac{m'd_2^2}{p^2}\Big).
\end{multline}
Using the Petersson trace formula, the $c'$-sum can be evaluated as a sum over $f\in B_k(l)$ plus a `diagonal' sum involving $\delta_{(n'=m)}$. The diagonal is easily estimated. For the sum over $f\in B_k(l)$, the Cauchy-Schwarz inequality and the spectral large sieve (see section \ref{specsieve}) then give that \eqref{remove} is bounded by $O(q^{-\half+\epsilon}p^{1+\epsilon})$, which is admissible. For details on the application of the large sieve, see its usage in the proof of Lemma \ref{bad-diag}.
\endproof 

Part of the proof of Lemma \ref{suffices} was relegated to the following result.

\begin{lemma}\label{bad-diag}
For any $d_1,d_2\le (qp)^{100}$, we have
\begin{align}
\label{bad-diag-main} q^\frac12   \sum_{\substack{n,m \ge 1\\ (m,q)=1}}  \frac{ \chi(n) \tau_{\overline{\chi}^2}(n)  }{n^\half}\frac{  \chi(m) \tau(m)   }{m^\half}V\Big(\frac{nd_1^2}{qp^2}\Big)V\Big(\frac{md_2^2}{qp^2}\Big) \sum_{g\in B_k(pq)} \frac{ \rho_g(nq)\overline{\rho_g(m)}}{\langle g, g\rangle }\ll q^\epsilon p^{1+\epsilon}.
 \end{align}
\end{lemma}
\proof
 Writing $B_k(pq)$ in terms of newforms using Lemma \ref{specialbasis}, we get that the left hand side of \eqref{bad-diag-main} equals
  \begin{align}
\label{newf} q^\frac12   \sum_{\substack{n,m \ge 1\\ (m,q)=1}}  \sum_{\substack{\delta \ell|pq}} \sum_{\substack{r_1|(\delta, nq) \\ r_2|(\delta, m)}}  \sum_{g\in B^*_k(\ell)}\frac{\nu_\delta(r_1)\nu_\delta(r_2) r_1^\half r_2^\half }{\langle g, g\rangle_{pq}}   \frac{ \chi(n) \tau_{\overline{\chi}^2}(n)  }{n^\half}\frac{  \chi(m) \tau(m)   }{m^\half} V\Big(\frac{n d_1^2}{qp^2}\Big)V\Big(\frac{m d_2^2}{qp^2}\Big)\rho_g\Big(\frac{nq}{r_1}\Big)\overline{\rho_g\Big(\frac{m}{r_2}\Big) },
 \end{align}
 for some coefficients $\nu_\delta(r_1), \nu_\delta(r_2)\ll (qp)^\epsilon$, where $\langle g, g\rangle_{pq}$ is the Petersson norm with respect to $\Gamma_0(pq)$. This normalization by $\langle g, g\rangle_{pq}$ serves to $L^2$-normalize with respect to $\Gamma_0(pq)$, so that we get an orthonormal basis as indicated in Lemma \ref{specialbasis}. By \eqref{obs}, we have $\langle g, g\rangle_{pq}\gg (qp)^{1-\epsilon}$.
 
 Consider first the case $q|\ell$ in \eqref{newf}. In this case we may replace $n$ by $nr_1$ (because $q\nmid \delta$, so $q\nmid r_1$, so the condition $r_1|nq$ implies $r_1|n$) and replace $m$ by $mr_2$ so that \eqref{newf} is bounded by
  \begin{align}
\label{newf-2} q^{-\half+\epsilon} p^{-1 +\epsilon}   \sum_{\substack{\delta \ell|pq\\ q|\ell}} \sum_{\substack{r_1|\delta \\ r_2|\delta }}  \sum_{g\in B^*_k(\ell)}   \Big|  \sum_{\substack{n,m \ge 1\\ (m,q)=1}} \frac{ \chi(nr_1) \tau_{\overline{\chi}^2}(nr_1)  }{n^\half}\frac{  \chi(mr_2) \tau(mr_2)   }{m^\half} V\Big(\frac{n r_1 d_1^2}{qp^2}\Big)V\Big(\frac{m r_2 d_2^2}{qp^2}\Big)\rho_g(nq)\overline{\rho_g(m)}\Big|.
 \end{align}
 We write $\rho_g(nq)=\rho_g(n)\rho_g(q)$ using Hecke multiplicativity and note that $\rho_g (q)\ll q^{-\half}$. The sum \eqref{newf-2} can then be estimated using the Cauchy-Schwarz inequality, the spectral large sieve, and \eqref{obs} to get the bound 
\begin{align*}
&q^{-1+\epsilon} p^{-1 +\epsilon}  \max_{\ell| pq} \Big( \sum_{g\in B^*_k(\ell)} \Big|\sum_{n\le q^{1+\epsilon}p^{2+\epsilon}} \frac{\alpha_n}{n^\half}\rho_g(n)\Big|^2\Big)^{\half} \Big( \sum_{g\in B^*_k(\ell)} \Big|\sum_{m\le q^{1+\epsilon}p^{2+\epsilon}} \frac{\beta_m}{m^\half}\rho_g(m)\Big|^2\Big)^{\half}\\
&\ll q^{-1+\epsilon} p^{-1 +\epsilon} (qp+qp^2)^\half(qp+qp^2)^\half\ll q^\epsilon p^{1+\epsilon},
\end{align*}
 for some complex numbers $\alpha_n, \beta_m$ bounded by $(qp)^\epsilon$.
 
 Now consider the case $q\nmid \ell$ in \eqref{newf}. Replacing $n$ by $n\frac{r_1}{(q,r_1)}$ and $m$ by $mr_2$, we get get that \eqref{newf} is bounded by
    \begin{multline}
\label{newf-3}  q^{-\half+\epsilon}p^{-1 +\epsilon}   \sum_{\substack{\delta \ell|pq\\ q\nmid \ell}} \sum_{\substack{r_1|\delta \\ r_2|\delta }}  \sum_{g\in B^*_k(\ell)}    \Big|  \sum_{\substack{n,m \ge 1\\ (m,q)=1}}   \frac{ (r_1,q)^{\half} \chi(\frac{nr_1}{(q,r_1)}) \tau_{\overline{\chi}^2}(\frac{nr_1}{(q,r_1)})  }{n^\half}\frac{  \chi(mr_2) \tau(mr_2)   }{m^\half} \\
V\Big(\frac{n r_1 d_1^2}{(q,r_1) qp^2 }\Big)V\Big(\frac{m r_2 d_2^2}{qp^2}\Big)\rho_g\Big(\frac{nq}{(r_1,q)}\Big)\overline{\rho_g(m)}\Big|.
 \end{multline}
 This sum is also ultimately treated using the spectral large sieve. However we do not apply it right away because we no longer have the saving $\rho_g(q)\ll q^{-\half}$ available as the level of $g$ is coprime to $q$. We only show the details for the most difficult case: the contribution to \eqref{newf-3} of the terms with $\ell=p, \delta=q, r_1=q, r_2=1$. This is the most difficult because, keeping in mind that the large sieve will give a bound in terms of the level and sum lengths, $\ell=p$ is the largest possible level under the condition $q\nmid \ell$, any choice other than $r_2=1$ would result in a shorter $m$-sum, $r_1=q$ introduces a large factor $(r_1,q)^\half$, and any choice other than $r_1=q$ would result in a shorter $n$-sum. The contribution of such terms is bounded by
    \begin{multline}
\label{newf-4}  q^{\epsilon}p^{-1 +\epsilon}   \sum_{g\in B^*_k(p)}    \Big|  \sum_{\substack{n,m \ge 1\\ (m,q)=1}}   \frac{ \chi(n) \tau_{\overline{\chi}^2}(n)  }{n^\half}\frac{  \chi(m) \tau(m)   }{m^\half} 
V\Big(\frac{n  d_1^2}{ qp^2 }\Big)V\Big(\frac{m d_2^2}{qp^2}\Big)\rho_g(n)\overline{\rho_g(m)}\Big|.
 \end{multline}
To treat this, we first remove the condition $(m,q)=1$ up to an admissible error. The contribution of the terms with $m=m'q$ is $O(q^\epsilon p^{1+\epsilon})$, using Hecke multiplicativity for $\rho_g(m'q)$, the Cauchy-Schwarz inequality, and the spectral large sieve. Next we use the definition of $V$ given in Lemma \ref{afe-use} to write
\begin{align}
\label{use-v} \sum_{m\ge 1} \frac{  \chi(m) \tau(m)   }{m^\half} V\Big(\frac{m d_2^2}{qp^2}\Big)\overline{\rho_g(m)} =  \int_{(2)}\int_{(2)}  \sum_{m_1, m_2 \ge 1} \frac{  \chi(m_1)\chi(m_2)  \overline{\rho_g(m_1m_2)}  }{m_1^{\half+s_1}m_2^{\half+s_2}}  \Big(\frac{qp^2}{d_2^2} \Big)^{s_1+s_2} G(s_1)G(s_2) ds_1 ds_2,
\end{align}
for some function $G(s)$ which is analytic for $\Re(s)>0$ and decays exponentially as $|\Im(s)|\to\infty$.
Using the Hecke multiplicativity relations, we get that \eqref{use-v} equals
\begin{align*}
&\int_{(2)}\int_{(2)}  \sum_{m_1, m_2, e \ge 1} \frac{ \mu(e)\chi^2(e)  \chi(m_1)\chi(m_2)  \overline{\rho_g(m_1)} \overline{\rho_g(m_2)}  }{e^{1+s_1+s_2} m_1^{\half+s_1}m_2^{\half+s_2}}  \Big(\frac{qp^2}{d_2^2} \Big)^{s_1+s_2}  G(s_1)G(s_2) ds_1 ds_2\\
&=\int_{(2)}\int_{(2)}   L(\thalf+s_1, \overline{g}\times \chi)L(\thalf+s_2, \overline{g}\times \chi) \sum_{ e \ge 1} \frac{ \mu(e)\chi^2(e)  }{e^{1+s_1+s_2} }  \Big(\frac{qp^2}{d_2^2} \Big)^{s_1+s_2}  G(s_1)G(s_2) ds_1 ds_2.
\end{align*}
The integrals can now be moved to the lines $\Re(s_i)=\epsilon$ and restricted to $|\Im(s_i)|\le (qp)^\epsilon$ by the exponential decay of $G(s)$. Similarly, the $n$-sum in \eqref{newf-4} can be written in terms of integrals involving $L(\thalf+s_1, g \times \chi)L(\thalf+s_2, g \times \overline{\chi})$ for $|s_i|\le (qp)^\epsilon$. Thus using Cauchy-Schwarz, we get that \eqref{newf-4} without the condition $(m,q)=1$, is bounded by
\[
\max_{|s|\le (qp)^\epsilon} q^{\epsilon}p^{-1 +\epsilon}   \sum_{g\in B^*_k(p)} |L(\thalf+s, g\times \chi)|^4.
\]
Since the analytic conductor of $L(\thalf+s, g\times \chi)$ is bounded by $q^\epsilon p^{2+\epsilon}$, by using a standard approximate functional equation, $|L(\thalf+s, g\times \chi)|^2$ can be expressed as a Dirichlet series of length $q^\epsilon p^{2+\epsilon}$. The spectral large sieve then gives a bound of
\[
q^{\epsilon}p^{-1 +\epsilon} (p+p^2)\ll q^\epsilon p^{1+\epsilon}.
\]
\endproof

{\it Remark.} Lemma \ref{bad-diag} gives rise to the error term $ O(q^{1+\epsilon}p^{1+\epsilon})$ in Theorem \ref{4th}. The main input in this lemma is the large sieve, and as such is not very deep. Thus with more effort, it may be possible to do better, at least in some ranges of $q$ and $p$. However we are content with this result because the other error term $O(p^{2+\epsilon} q^{\half+\vartheta+\epsilon})$ in Theorem \ref{4th} is anyway also $O(q^{1+\epsilon}p^{1+\epsilon})$ on the Ramanujan Conjecture when $p^{2-\epsilon}<q< p^{2+\epsilon}$, which is our main range of interest.

\

Now we return to Lemma \ref{suffices}, and restrict the sum to dyadic intervals. It suffices to prove that 

\begin{multline*}
 \frac{1}{N^{\half}M^{\half} Cq^{\frac12}} \sum_{\substack{n,m,c \ge 1\\ (c,pq)=1}}  \chi(n) \tau_{\overline{\chi}^2}(n)   \chi(m) \tau(m)  S(n,m\overline{q},c) J_{k-1}\Big(\frac{4\pi \sqrt{nm}}{c\sqrt{q}}\Big)  V_1\Big(\frac{n}{N}\Big)  V_2\Big(\frac{m}{M}\Big)  V_3\Big(\frac{c}{C}\Big)\\
 \ll (qp)^\epsilon p^2q^{-\half+\vartheta},
\end{multline*}
for any 
\begin{align}
\label{nmranges} 1\le N,M\le q^{1+\epsilon}p^{2+\epsilon}
\end{align}
and $C \ge 1$, and weight functions $V_i$ compactly supported on $(1,2)$ and satisying $V_i^{(j)}(x)\ll_{j} (qp)^\epsilon$. Note that these weight functions are arbitrary, and we will redefine them (and others like them) from line to line to suit our needs. For example, $\sqrt{n} V(\frac{n}{N})$ equals $\sqrt{N} \sqrt{\frac{n}{N}} V(\frac{n}{N})$, which may simply be written as $\sqrt{N}V(\frac{n}{N})$ after redefining $\sqrt{x}V(x)$ to equal $V(x)$.

Using Lemma \ref{jlemma}, we need to establish
\begin{multline*}
\frac{1}{C^2q} \sum_{\substack{n,m,c \ge 1\\ (c,pq)=1}}  \chi(n) \tau_{\overline{\chi}^2}(n)   \chi(m) \tau(m)  S(n,m\overline{q},c) W_1\Big(\frac{\sqrt{nm}}{c\sqrt{q}}\Big)  V_1\Big(\frac{n}{N}\Big)  V_2\Big(\frac{m}{M}\Big)  V_3\Big(\frac{c}{C}\Big)
 \ll (qp)^\epsilon p^2q^{-\half+\vartheta}
\end{multline*}
in the case that $\frac{\sqrt{NM}}{C\sqrt{q}}\le (qp)^\epsilon$, and 
\begin{multline*}
 \frac{1}{N^\frac34 M^\frac34 C^\half q^\frac14}\sum_{\substack{n,m,c \ge 1\\ (c,pq)=1}}  \chi(n) \tau_{\overline{\chi}^2}(n)   \chi(m) \tau(m)  S(n,m\overline{q},c) e\Big(\frac{\pm 2\sqrt{nm}}{c\sqrt{q}}\Big) W_2\Big(\frac{ \sqrt{nm}}{c\sqrt{q}}\Big)  V_1\Big(\frac{n}{N}\Big)  V_2\Big(\frac{m}{M}\Big)  V_3\Big(\frac{c}{C}\Big)\\
 \ll (qp)^\epsilon p^2q^{-\half+\vartheta}
\end{multline*}
in the case that $\frac{\sqrt{NM}}{C\sqrt{q}}> (qp)^\epsilon$, some some real functions $W_1$ and $W_2$ satisfying bounds as given in the lemma. (We are assuming these are real by taking real and imaginary parts; this is also why the $\pm$ arises in the exponential factor.) We may write $W_1(\frac{\sqrt{nm}}{c\sqrt{q}})$ and $W_2(\frac{\sqrt{nm}}{c\sqrt{q}})$ in terms of their real and imaginary parts, use Lemma \ref{mellin-split} to separate the variables $n,m$, and $c$, and then absorb the resulting factors $n^{-\frac{s}{2}}, m^{-\frac{s}{2}}$, and $c^{s}$ into the functions $V_1, V_2$, and $V_3$. Thus it suffices to prove
\begin{align}
\label{nonosc} \frac{1}{C^2q} \sum_{\substack{n,m,c \ge 1\\ (c,pq)=1}}  \chi(n) \tau_{\overline{\chi}^2}(n)   \chi(m) \tau(m)  S(n,m\overline{q},c)   V_1\Big(\frac{n}{N}\Big)  V_2\Big(\frac{m}{M}\Big)  V_3\Big(\frac{c}{C}\Big)
 \ll (qp)^\epsilon p^2q^{-\half+\vartheta}
\end{align}
for any
$
N\le q^{1+\epsilon}p^{2+\epsilon}, \ \ M\le q^{1+\epsilon}p^{2+\epsilon}, \ \  \frac{\sqrt{NM}}{C\sqrt{q}}\le (qp)^\epsilon,
$
and
\begin{multline}
 \label{osc}  \frac{1}{N^\frac34 M^\frac34 C^\half q^\frac14} \sum_{\substack{n,m,c \ge 1\\ (c,pq)=1}}  \chi(n) \tau_{\overline{\chi}^2}(n)   \chi(m) \tau(m)  S(n,m\overline{q},c) e\Big(\frac{\pm 2\sqrt{nm}}{c\sqrt{q}}\Big)  V_1\Big(\frac{n}{N}\Big)  V_2\Big(\frac{m}{M}\Big)  V_3\Big(\frac{c}{C}\Big)\\
 \ll (qp)^\epsilon p^2q^{-\half+\vartheta}
\end{multline}
for any
$
N\le q^{1+\epsilon}p^{2+\epsilon}, \ \ M\le q^{1+\epsilon}p^{2+\epsilon}, \ \  \frac{\sqrt{NM}}{C\sqrt{q}}> (qp)^\epsilon.
$
We may also assume the crude bound
\begin{align}
\label{crude} C\le (qp)^{50}.
\end{align}
This is clear in the case $\frac{\sqrt{NM}}{C\sqrt{q}}> (qp)^\epsilon$, while in the complementary case we may make this assumption by summing trivially using the Weil bound in \eqref{nonosc}.


\section{The oscillatory case}
In this section we prove \eqref{osc}. This case is more involved than \eqref{nonosc} analytically because the weight function in the sum has the oscillatory factor $e(\frac{2\sqrt{nm}}{c\sqrt{q}})$.
We keep in mind that in this case, we have
\begin{align}
\label{ccase2} \frac{\sqrt{NM}}{C\sqrt{q}}> (qp)^\epsilon.
\end{align}


\subsection{Application of Voronoi summation} \label{appvor1}

We write
\[
S(n,m\overline{q},c)=\summ_{b\bmod c} e\Big(\frac{nb+m\overline{qb}}{c}\Big),
\]
and apply Lemma \ref{vor} to the $n$ and $m$ sums in \eqref{osc}. Thus we apply Voronoi summation twice, and after a substitution in the integral transforms, we get that it suffices to prove
\begin{multline}
 \label{need4}    \frac{N^\frac14 M^\frac14}{ C^\half q^\frac14} \sum_{\substack{n,m,c \ge 1\\ (c,pq)=1}} \frac{\chi(c^2)}{c^2 p^2} \overline{\chi}(n) \tau_{\overline{\chi}^2}(n)  \chi(m) \tau(m)  S(\pm nq,m\overline{p}^4,c) \\
  \int_0^\infty \int_0^\infty e\Big( \frac{\pm 2\sqrt{NMxy}}{c\sqrt{q}}\Big) B_1\Big(\frac{4\pi\sqrt{nNx}}{cp}\Big)B_2\Big(\frac{4\pi\sqrt{mMy}}{cp}\Big) V_1(x)V_2(y) V_3\Big( \frac{ c}{C}\Big) dxdy
 \ll (qp)^\epsilon p^2q^{-\half+\vartheta},
\end{multline}
for each combination of Bessel functions $\{B_1, B_2\}=\{ Y_0, K_0\}$ and for each choice of sign $\pm$. The two occurences of $\pm$ are independent of each other, but when $B_1=B_2=Y_0$, then the sign inside the Kloosterman sum is positive.

 We claim that we need only consider the case $B_1=B_2=Y_0$ in \eqref{need4}. For if $B_1=K_0$, say, then either $\frac{\sqrt{nN}}{Cp}\ge (qp)^\epsilon$, in which case $K_0$ is tiny by Lemma \ref{jlem} and the bound in \eqref{need4} follows, or $\frac{\sqrt{nN}}{Cp}< (qp)^\epsilon$. But in this case we can integrate by parts repeatedly in \eqref{need4} with respect to $x$, integrating $e( \frac{\pm 2\sqrt{NMxy}}{c\sqrt{q}})$ and differentiating the rest. Using Lemma \ref{jlem} to control the size of derivatives of $K_0$, this will show that the integral is bounded by $(qp)^{-100}$ so that \eqref{need4} is established. Thus we can assume $B_1=B_2=Y_0$, in which case the sign inside the Kloosterman sum in \eqref{need4} is $+$. We can also assume that  
\[
\frac{\sqrt{nN}}{cp}\ge (qp)^\epsilon, \ \ \ \ \ \frac{\sqrt{mM}}{cp}\ge (qp)^\epsilon,
\]
for if $\frac{\sqrt{nN}}{cp}<(qp)^\epsilon$, say, then we can once again integrate by parts, using Lemma \ref{jlem} to control the size of derivatives of $Y_0$, to show the bound required in \eqref{need4}. Now, using \eqref{y0long}, we need to prove 
\begin{multline*}
 \frac{N^\frac14 M^\frac14}{ C^\half q^\frac14} \sum_{\substack{n,m,c \ge 1\\ (c,pq)=1}} \frac{\chi(c^2)}{c^2 p^2} \Big(\frac{\sqrt{mM}}{cp}\frac{\sqrt{nN}}{cp}\Big)^{-\half}  \overline{\chi}(n) \tau_{\overline{\chi}^2}(n)  \chi(m) \tau(m)  S( nq,m\overline{p}^4,c) 
  \int_0^\infty \int_0^\infty e\Big( \frac{\pm 2\sqrt{NMxy}}{c\sqrt{q}}\Big) \\
 e\Big(\frac{\pm 2\sqrt{nNx}}{cp}\Big)e\Big(\frac{\pm 2\sqrt{mMy}}{cp}\Big)
  W_1\Big(\frac{2\sqrt{nNx}}{cp}\Big)W_2\Big(\frac{2\sqrt{mMy}}{cp}\Big) V_1(x)V_2(y) V_3\Big( \frac{c}{C}\Big) dxdy
 \ll (qp)^\epsilon p^2q^{-\half+\vartheta},
\end{multline*}
for some real functions $W_1, W_2$ satisfying $W_i^{(j)}(x)\ll_j x^{-j}$ for $j\ge 0$, and where each sign $\pm$ is independent. Making the substituting $\sqrt{x}\to x$ and $\sqrt{y}\to y$, and redefining the weight functions, we need to show
\begin{multline}
\label{need5} \frac{p}{q^\frac14 C^\frac12 } \sum_{\substack{n,m,c \ge 1\\ (c,pq)=1}}\frac{1}{n^\frac14 m^\frac14} \overline{\chi}(n) \tau_{\overline{\chi}^2}(n)  \chi(m) \tau(m)  \frac{\chi(c^2)S( nq,m\overline{p}^4,c)}{cp^2} 
  \int_0^\infty \int_0^\infty e\Big( \frac{\pm 2\sqrt{NM}xy}{c\sqrt{q}}\Big) \\
 e\Big(\frac{\pm 2\sqrt{nN}x}{cp}\Big)e\Big(\frac{\pm 2\sqrt{mM}y}{cp}\Big)
  W_1\Big(\frac{2\sqrt{nN}x}{cp}\Big)W_2\Big(\frac{2\sqrt{mM}y}{cp}\Big) V_1(x)V_2(y) V_3\Big( \frac{c}{C}\Big) dxdy
 \ll (qp)^\epsilon p^2q^{-\half+\vartheta}.
\end{multline}
Writing $Z_1(x)= W_1(\frac{2\sqrt{nN}x}{cp}) V_1(x)$, $Z_2(y)= W_2(\frac{2\sqrt{nN}y}{cp}) V_2(y)$ and evaluating the $x$-integral in terms of the Fourier transform $\hat{Z_1}$, we need to show
\begin{multline*}
 \frac{p}{q^\frac14 C^\frac12 } \sum_{\substack{n,m,c \ge 1\\ (c,pq)=1}}\frac{1}{n^\frac14 m^\frac14}\overline{\chi}(n) \tau_{\overline{\chi}^2}(n)  \chi(m) \tau(m) \frac{\chi(c^2)S( nq,m\overline{p}^4,c)}{cp^2} 
 \\
 \int_0^\infty  e\Big(\frac{\pm 2\sqrt{mM}y}{cp}\Big) \hat{Z}_1\Big(\frac{\pm 2\sqrt{NM} y}{c\sqrt{q}} \pm  \frac{2\sqrt{nN} }{cp} \Big)
 Z_2(y) V_3\Big( \frac{c}{C}\Big) dy
 \ll (qp)^\epsilon p^2q^{-\half+\vartheta}.
\end{multline*}
Now we make the substitution $\frac{\pm 2\sqrt{NM} y}{c\sqrt{q}} \pm  \frac{2\sqrt{nN} }{cp}\to y$ to see that what we need to prove is
\begin{multline}
\label{need6} \frac{q^\frac14 p C^\half }{N^\frac12 M^\frac12  } \sum_{\substack{n,m,c \ge 1\\ (c,pq)=1}}\frac{1}{n^\frac14 m^\frac14}  \overline{\chi}(n) \tau_{\chi^2}(n)  \chi(m) \tau(m)  \frac{\chi(c^2)S( nq,m\overline{p}^4,c)}{cp^2} 
 \\
e\Big(\frac{\pm 2\sqrt{nmq}}{cp^2}\Big)   \int_0^\infty  e\Big(\frac{\pm \sqrt{mq}}{\sqrt{N} p}y\Big) \hat{Z_1}(y)
 Z_2\Big(\frac{\pm c\sqrt{q}y}{2\sqrt{NM}} \pm \frac{\sqrt{nq}}{p\sqrt{M}}\Big) V_3\Big( \frac{c}{C}\Big) dy
 \ll (qp)^\epsilon p^2q^{-\half+\vartheta}.
\end{multline}
By repeated integration by parts, we have
\begin{align}
\label{fourierbound} \hat{Z_1}^{(j)}(y)= \frac{d^j}{dy^j} \intt e(-xy)Z_1(x) dx\ll_{j,l} (qp)^\epsilon \min(1, y^{-l})
\end{align}
for any $j, l\ge 0$, so that that we may restrict the integral to $|y|\ll (qp)^\epsilon$. 

We end this section by observing that in \eqref{need6}, we may restrict to
\begin{align}
\label{xrange} (qp)^{-\epsilon} \frac{\sqrt{NM}}{C\sqrt{q}} \le \frac{ \sqrt{nmq}}{cp^2}\le  (qp)^\epsilon \frac{\sqrt{NM}}{C\sqrt{q}}.
\end{align}
This is because by repeatedly integrating by parts the $y$-integral in \eqref{need6}, integrating the exponential and differentiating the rest, we can assume that
\[
\frac{\sqrt{mq}}{\sqrt{N} p} \le (qp)^\epsilon,
\]
otherwise the left hand side of \eqref{need6} is bounded by $(qp)^{-100}$. This uses \eqref{ccase2}. Similarly, by returning to \eqref{need5} and evaluating the $y$-integral instead of the $x$-integral, we can assume that 
\[
\frac{\sqrt{nq}}{\sqrt{M} p} \le (qp)^\epsilon.
\]
Next we show that we can assume the lower bounds 
\[
\frac{\sqrt{mq}}{\sqrt{N} p} \ge (qp)^{-\epsilon}, \ \ \ \frac{\sqrt{nq}}{\sqrt{M} p} \ge (qp)^{-\epsilon},
\]
after which \eqref{xrange} will follow. Suppose that these lower bounds don't hold; say 
\[
\frac{\sqrt{nq}}{\sqrt{M} p} < (qp)^{-\epsilon}.
\]
Then by the support of $Z_2$, we have that $ \frac{\sqrt{NM}}{C\sqrt{q}}\ll y\ll  \frac{\sqrt{NM}}{C\sqrt{q}}$ in  in \eqref{need6}. But then by \eqref{ccase2} and \eqref{fourierbound}, the left hand side of \eqref{need6} is bounded by $(qp)^{-100}$.


\subsection{Application of Kuznetsov's formula}\label{appkuz}

By \eqref{xrange}, we can insert in in \eqref{need6} a function $\psi(\frac{4\pi \sqrt{nmq}}{cp^2})$, such that $\psi(x)$ is smooth, supported on $X<x<2X$, and satisfies $\psi^{(j)}(x)\ll X^{-j}$, with
\[
(qp)^{-\epsilon} \frac{\sqrt{NM}}{C\sqrt{q}} \le X \le (qp)^\epsilon \frac{\sqrt{NM}}{C\sqrt{q}}.
\]
Thus in \eqref{need6}, we may replace $(nm)^{-\frac14}$ by its typical size $(N M)^{-\frac14} p^{-1} q^{\half}$, by redefining $\psi$. We may apply Lemma \ref{mellin-split} to  and $Z_2(x \pm \frac{\sqrt{nq}}{p\sqrt{M}})$ and $V_3(x)$ to separate factors $c^{-s}$ for $|s|<(qp)^\epsilon$ from $ Z_2(\frac{\pm c\sqrt{q}y}{2\sqrt{NM}} \pm \frac{\sqrt{nq}}{p\sqrt{M}}) V_3( \frac{c}{C})$, and absorb them into $\psi(\frac{4\pi \sqrt{nmq}}{cp^2})$. Thus it suffices to prove
\begin{multline}
\label{need7a} \frac{q^\frac34 C^\half }{N^\frac34 M^\frac34}  \sum_{\substack{n\asymp \frac{M p^2}{q}\\ m\asymp \frac{N p^2}{q}}} \sum_{\substack{c \ge 1\\ (c,p)=1}}\alpha_n \beta_m \frac{\chi(c^2)S( nq,m\overline{p}^4,c)}{cp^2}
  e\Big(\frac{\pm 2\sqrt{nmq}}{cp^2}\Big) \psi\Big(\frac{4\pi \sqrt{nmq}}{cp^2}\Big)  \ll (qp)^\epsilon p^2q^{-\half+\vartheta},
 \end{multline}
 and
 \begin{multline}
 \label{need7b}  \frac{q^\frac34 C^\half }{N^\frac34 M^\frac34 }  \sum_{\substack{n\asymp \frac{M p^2}{q}\\ m\asymp \frac{N p^2}{q}}} \sum_{\substack{c \ge 1\\ (c,p)=1}}\alpha_n \beta_m   \frac{\chi(c^2)S( nq,m\overline{p}^4,cq)}{cqp^2} 
 e\Big(\frac{\pm 2\sqrt{nmq}}{cqp^2}\Big) \psi\Big(\frac{4\pi \sqrt{nmq}}{cqp^2}\Big)  \ll (qp)^\epsilon p^2q^{-\half+\vartheta},
 \end{multline}
 where $\alpha_n, \beta_m$ are arbitrary complex numbers bounded by $(qp)^\epsilon$ and where we write $n\asymp \frac{M p^2}{q}$ to mean $(qp)^{-\epsilon} \frac{M p^2}{q}\le n\le (qp)^\epsilon \frac{M p^2}{q}$, and so on (for these inequalities, see the end of section \ref{appvor1}). We have two sums because we used inclusion-exclusion to remove the condition $(c,q)=1$ from \eqref{need6}.
 
We insert $e(\frac{m}{p^2})$ in \eqref{need7a} and $e(\frac{m\overline{q}}{p^2})$ in \eqref{need7b}, which we can do since $\beta_m$ are anyway arbitrary, and then apply Theorem \ref{kuznetsov}. Thus we are reduced to showing
 \begin{align}
\label{need8a}  \frac{q^\frac34 C^\half }{N^\frac34 M^\frac34}  \Big( \sum_{\substack{g \in B_0(p^4,\chi^2)}}  \frac{\tilde{\phi}(t_g) }{\cosh(\pi t_g)}  \sum_{\substack{n\asymp \frac{M p^2}{q}\\ m\asymp \frac{N p^2}{q}}} \frac{\rho_{g,\infty}(nq)  \overline{\rho_{g,1}(m)}}{\| g\|^2} +\mathcal{H}+\mathcal{E}\Big)
  \ll (qp)^\epsilon p^2q^{-\half+\vartheta},
 \end{align}
 and
  \begin{align}
\label{need8b}  \frac{q^\frac34 C^\half }{N^\frac34 M^\frac34}  \Big( \sum_{\substack{g \in B_0(qp^4,\chi^2)}}  \frac{\tilde{\phi}(t_g) }{\cosh(\pi t_g)}  \sum_{\substack{n\asymp \frac{M p^2}{q}\\ m\asymp \frac{N p^2}{q}}} \frac{\rho_{g,\infty}(nq)  \overline{\rho_{g,1}(m)}}{\| g\|^2} +\mathcal{H}+\mathcal{E}\Big)
  \ll (qp)^\epsilon p^2q^{-\half+\vartheta},
 \end{align}
 where 
 \[
 \phi(x)=e^{\pm i x} \psi(x)
 \]
satisfies $\phi^{(j)}(x)\ll 1$, and we have only written out the Maass form contribution $\mathcal{M}$. We will not show the details for $\mathcal{H}$ and $\mathcal{E}$ because their treatment is similar and easier because there are no issues with the Ramanujan conjecture.
 
 We only show the details for the Maass form contribution in \eqref{need8a}, as  \eqref{need8b} is similarly treated. By Lemma \ref{transf-bounds} (c), for real $t_g$ we can restrict to 
\[
t_g\ll (qp)^\epsilon \Big(\frac{\sqrt{NM}}{C\sqrt{q}}\Big)^{\half},
\]
in which range we have
\[
\tilde{\phi}(t_g)\ll (qp)^\epsilon \Big(\frac{\sqrt{NM}}{C\sqrt{q}}\Big)^{-\half}.
\]
Using the Cauchy-Schwarz inequality, for the Maass contribution of real $t_f$ to \eqref{need8a} we need to show
\begin{multline}
\label{need9a} \frac{qC}{NM}\Big(\sum_{\substack{g\in B_0(p^4,\chi^2)\\ t_g \le (qp)^\epsilon (\frac{\sqrt{NM}}{C\sqrt{q}})^{\half} }}  \frac{1}{\| g \|^2} \frac{1}{\cosh(\pi t_g)}\Big| \sum_{n\asymp \frac{ M p^2}{q}} \alpha_n \rho_{g,\infty}(nq)\Big|^2\Big)^\half \\
\times\Big( \sum_{\substack{g\in B_0(p^4,\chi^2)\\ t_g \le (qp)^\epsilon (\frac{\sqrt{NM}}{C\sqrt{q}})^{\half} }} \frac{1}{\| g \|^2} \frac{1}{\cosh(\pi t_g)}\Big| \sum_{m\asymp \frac{ N p^2}{q}} \overline{\beta_m}   \rho_{g,1}(m) \Big|^2\Big)^\half  \ll (qp)^\epsilon p^2q^{-\half+\vartheta}.
\end{multline}
 We first see what we should expect: if we were working with newforms and had $ \rho_{g,\infty}(n) \rho_{g,\infty}(q) $ instead of $\rho_{g,\infty}(nq)$, then we could use $ \rho_{g,\infty}(q)\ll q^{\vartheta+\epsilon}$ and the spectral large sieve to get that the left hand side of \eqref{need9a} is bounded by
\begin{align}
\label{largesieveexpectation} & q^{\vartheta+\epsilon} p^\epsilon \frac{qC}{NM} \Bigg( \Big(\frac{\sqrt{NM}}{C\sqrt{q}} +\frac{1}{p^4} \frac{Mp^2}{q}\Big) \frac{Mp^2}{q}  \Bigg)^\half  \Bigg( \Big(\frac{\sqrt{NM}}{C\sqrt{q}}  + \frac{1}{p^4}\frac{Np^2}{q}\Big) \frac{Np^2}{q} \Bigg)^\half .
\end{align}
We have that $\max( \frac{1}{p^4} \frac{Mp^2}{q}, \frac{1}{p^4} \frac{Np^2}{q}) \le (qp)^\epsilon \frac{\sqrt{NM}}{C\sqrt{q}} $ because $\max(N,M)\le q^{1+\epsilon}p^{2+\epsilon}$ by \eqref{nmranges} and $\frac{\sqrt{NM}}{C\sqrt{q}}\ge (qp)^\epsilon$ by \eqref{ccase2}. Thus we get that is $O(p^{2+\epsilon}q^{-\half+\vartheta+\epsilon})$, as required.  In order to rigorously separate $n$ and $q$, we write $B_0(p^4,\chi^2)$ in terms of newforms using Lemma \ref{specialbasis}, to get that
\begin{align*}
&\sum_{\substack{g\in B_0(p^4,\chi^2)\\ t_g \le (qp)^\epsilon (\frac{\sqrt{NM}}{C\sqrt{q}})^{\half} }}  \frac{1}{\| g \|^2} \frac{1}{\cosh(\pi t_g)}\Big| \sum_{n\asymp \frac{ M p^2}{q}} \alpha_n \rho_{g,\infty}(nq)\Big|^2\\
&=\sum_{\substack{\delta\ell|p^4}} \sum_{\substack{g\in B_0^*(\ell,\chi^2)\\ t_g \le (qp)^\epsilon (\frac{\sqrt{NM}}{C\sqrt{q}})^{\half} }} \frac{1}{\langle g, g\rangle_{p^4}^2 } \frac{1}{\cosh(\pi t_g)}\Big| \sum_{n\asymp \frac{ M p^2}{q}} \sum_{r| (nq,\delta)} \nu_\delta(r) r^\half \alpha_n \rho_{g,\infty}\Big(\frac{nq}{r}\Big)\Big|^2
\end{align*}
for some coefficients $\nu_\delta(r)\ll (qp)^\epsilon$, where $\langle g, g\rangle_{p^4}$ is the Petersson norm with respect to $\Gamma_0(p^4)$. We have $ \langle g, g\rangle_{\ell} (\frac{p^4}{\ell})^{-\epsilon} \ll \langle g, g\rangle_{p^4} \ll \langle g, g\rangle_{\ell} (\frac{p^4}{\ell})^\epsilon$ by the remarks at the end of section \ref{specsieve}. We also have 
\begin{align*}
\Big| \sum_{n\asymp \frac{ M p^2}{q}} \sum_{r| (nq,\delta)} \nu_\delta(r) r^\half \alpha_n \rho_{g,\infty}\Big(\frac{nq}{r}\Big)\Big|^2 
\ll (qp)^\epsilon  \Big|  \sum_{r|\delta} r^{\half} \sum_{n\asymp \frac{ M p^2}{qr}} \alpha_{nr} \rho_{g,\infty}(nq) \Big|^2
\end{align*}
after replacing $n$ by $nr$, and 
\begin{align*}
 \sum_{n\asymp \frac{ M p^2}{qr}} \alpha_{nr} \rho_{g,\infty}(nq)= \sum_{d|(n,q)}  \sum_{n\asymp \frac{ M p^2}{qr}} \mu(d)\chi^2(d) \alpha_{nr} \rho_{g,\infty}\Big(\frac{n}{d}\Big) \rho_{g,\infty}\Big(\frac{q}{d}\Big) 
\ll q^{\vartheta+\epsilon} \sum_{d|q}  \Big|  \sum_{n\asymp \frac{ M p^2}{qdr}} \alpha_{ndr} \rho_{g,\infty}(n) \Big|
\end{align*}
using the Hecke relations. Thus we get that the square of the left hand side of \eqref{need9a} is bounded by
\begin{align*}
(qp)^\epsilon\Big(\frac{qC}{NM}\Big)^2q^{2\vartheta+\epsilon} \sum_{\substack{ \delta \ell | p^4  \\ r_1|\delta, r_2 | \delta \\ d_1|q, d_2|q}} &r_1^\half r_2^\half \frac{\ell}{p^4} \Big(\sum_{\substack{g\in B_0^*(\ell,\chi^2)\\ t_g \le (qp)^\epsilon (\frac{\sqrt{NM}}{C\sqrt{q}})^{\half} }}  \frac{1}{\| g \|^2} \frac{1}{\cosh(\pi t_g)}\Big| \sum_{n_1\asymp \frac{ M p^2}{q d_1 r_1 }} \alpha_{n_1d_1r_2} \rho_{g,\infty}(n_1)\Big)\Big|^2\Big)^\half \\
 &\times \Big(\sum_{\substack{g\in B_0^*(\ell,\chi^2)\\ t_g \le (qp)^\epsilon (\frac{\sqrt{NM}}{C\sqrt{q}})^{\half} }}  \frac{1}{\| g \|^2} \frac{1}{\cosh(\pi t_g)}\Big| \sum_{n_2\asymp \frac{ M p^2}{q d_2 r_2 }} \alpha_{n_2d_2r_2} \rho_{g,\infty}(n_2) \rho_{g,\infty}\Big|^2\Big)^\half \\
&\times\Big( \sum_{\substack{g\in B_0(p^4,\chi^2)\\ t_g \le (qp)^\epsilon (\frac{\sqrt{NM}}{C\sqrt{q}})^{\half} }} \frac{1}{\| g \|^2} \frac{1}{\cosh(\pi t_g)}\Big| \sum_{m\asymp \frac{ N p^2}{q}} \overline{\beta_m}   \rho_{g,1}(m) \Big|^2\Big).
\end{align*}
Now applying the large sieve gives the bound $O(p^{2+\epsilon}q^{-\half+\vartheta+\epsilon})$ as we expected.

Now we return to \eqref{need8a} and consider the contribution the exceptional eigenvalues. We have by  Lemma \ref{transf-bounds} (b) with $Z=\frac{\sqrt{NM}}{C\sqrt{q}}$ that $\tilde{\phi}(t_g)\ll (qp)^\epsilon$. Thus using the Cauchy-Schwarz inequality, the Maass contribution of imaginary $t_g$ in \eqref{need8a} is bounded by
\begin{align}
\label{imag} \frac{q^\frac34 C^\half }{N^\frac34 M^\frac34}  &\Big(\sum_{\substack{g\in B_0(p^4,\chi^2)\\ t_g \ll 1 }}  \frac{1}{\| g \|^2} \frac{1}{\cosh(\pi t_g)}\Big| \sum_{n\asymp \frac{ M p^2}{q}} \alpha_n \rho_{g,\infty}(nq)\Big|^2\Big)^\half \\
\nonumber  \times &\Big( \sum_{\substack{g\in B_0(p^4,\chi^2)\\ t_g \ll 1 }} \frac{1}{\| g \|^2} \frac{1}{\cosh(\pi t_g)}\Big| \sum_{m\asymp \frac{ N p^2}{q}} \overline{\beta_m}   \rho_{g,1}(m) \Big|^2\Big)^\half.
\end{align}
If instead of $\rho_{g,\infty}(nq)$ we had $\rho_{g,\infty}(n)\rho_{g,\infty}(q)$, then by the spectral large sieve and the Ramanujan bound for $\rho_g(q)$, we would get that \eqref{imag} is bounded by
\begin{align}
\label{imag2} q^{\vartheta+\epsilon} p^\epsilon \frac{q^\frac34 C^\half }{N^\frac34 M^\frac34}  \Bigg( \Big(1 +\frac{1}{p^4} \frac{Mp^2}{q}\Big) \frac{Mp^2}{q}  \Bigg)^\half  \Bigg( \Big(1 + \frac{1}{p^4}\frac{Np^2}{q}\Big) \frac{Np^2}{q} \Bigg)^\half .
\end{align}
Using again that $\max(N,M)\le q^{1+\epsilon}p^{2+\epsilon}$, we get that \eqref{imag2} is bounded by 
\begin{align*}
 q^{-\half+\vartheta+\epsilon} p^{2+\epsilon}  \Big(\frac{\sqrt{q}C}{\sqrt{NM}} \Big)^\half\ll q^{-\half+\vartheta+\epsilon} p^{2+\epsilon},
 \end{align*}
using \eqref{ccase2}. We can rigorously separate $n$ and $q$ just as we did for real $t_g$ and arrive at the same conclusion.


\section{The non-oscillatory case}
In this section we prove \eqref{nonosc}. Unlike the previous case, the weight function in the sum does not have an oscillatory factor, and thus is somewhat easier. We will keep the details concise.

Keep in mind that in this case, we have
\begin{align}
\label{ccase1} q^{-100} \le \frac{\sqrt{NM}}{C\sqrt{q}}\le (qp)^\epsilon,
\end{align}
using \eqref{crude}.
Also recall that
\begin{align}
\label{eq3} \frac{\sqrt{NM}}{qp^2}<(qp)^\epsilon
\end{align}
 using $N,M\le (qp)^\epsilon qp^2$.

\subsection{Application of Voronoi summation}
On applying the Voronoi summation formula twice in \eqref{nonosc} and making a substitution in the resulting integral transforms, we need to prove
\begin{multline}
 \label{need10}    \frac{NM}{ C^2 q} \sum_{\substack{n,m,c \ge 1\\ (c,pq)=1}} \frac{\chi(c^2)}{c^2 p^2} \overline{\chi}(n) \tau_{\overline{\chi}^2}(n)  \chi(m) \tau(m)  S(\pm nq,m\overline{p}^4,c) \\
  \int_0^\infty \int_0^\infty B_1\Big(\frac{4\pi\sqrt{nNx}}{cp}\Big)B_2\Big(\frac{4\pi\sqrt{mMy}}{cp}\Big) V_1(x)V_2(y) V_3\Big( \frac{ c}{C}\Big) dxdy
 \ll (qp)^\epsilon p^2q^{-\half+\vartheta},
\end{multline}
for each combination of Bessel functions $\{B_1, B_2\}=\{ Y_0, K_0\}$ and for each choice of sign $\pm$.  
We claim that we can assume
\begin{align*}
n\le (qp)^\epsilon \frac{C^2 p^2}{N}, \ \ \ \ m\le (qp)^\epsilon \frac{C^2 p^2}{M}.
\end{align*}
Clearly this follows for $n$ when $B_1=K_0$, by Lemma \ref{jlem}. If $B_1=Y_0$ and $n> (qp)^\epsilon \frac{(Cp)^2}{N}$, then we may apply \eqref{ybes-big} and integrate by parts repeatedly in \eqref{need4} with respect to $x$ to see that the contribution of such terms is less than $(qp)^{-100}$.  Arguing the same way for $m$ and $B_2$, the claim follows. Thus we can work in dyadic intervals of $n$ and $m$, restricting the range of summation in \eqref{need10} to $N'\le n\le 2N'$ and $M'\le m\le 2M'$, where
\begin{align}
\label{n'm'range} 1\le N'\le (qp)^\epsilon \frac{C^2p^2}{N}, \ \ 1\le M'\le (qp)^\epsilon \frac{C^2p^2}{M}.
\end{align}
Then we can insert a function $\phi(\frac{4\pi \sqrt{nmq}}{cp^2})$, such that $\phi(x)$ is smooth, supported on $X<x<2X$, and satisfies $\phi^{(j)}(x)\ll X^{-j}$, with
\[
 (qp)^{-\epsilon} \frac{\sqrt{N'M'q}}{Cp^2} \le X \le (qp)^\epsilon \frac{\sqrt{N'M'q}}{Cp^2}.
\]
Thus it suffices to prove
\begin{align}
\label{need11} \frac{NM}{C^3 q}  \sum_{\substack{N' \le n\le 2N' \\ M' \le n\le 2M' }} \sum_{\substack{c \ge 1\\ (c,p)=1}}\alpha_n \beta_m \frac{\chi(c^2)S( \pm nq,m\overline{p}^4,c)}{cp^2}
 \phi\Big(\frac{4\pi \sqrt{nmq}}{cp^2}\Big)  \ll (qp)^\epsilon p^2q^{-\half+\vartheta},
 \end{align}
 and
\begin{align*}
\frac{NM}{C^3 q}  \sum_{\substack{N' \le n\le 2N' \\ M' \le n\le 2M' }} \sum_{\substack{c \ge 1\\ (c,p)=1}}\alpha_n \beta_m \frac{\chi(c^2)S( \pm nq,m\overline{p}^4,cq)}{cqp^2}
 \phi\Big(\frac{4\pi \sqrt{nmq}}{cqp^2}\Big)  \ll (qp)^\epsilon p^2q^{-\half+\vartheta},
 \end{align*}
 where $\alpha_n, \beta_m$ are arbitrary complex numbers bounded by $(qp)^\epsilon$. We have two sums because we used inclusion-exclusion to remove the condition $(c,q)=1$ from \eqref{need6}.
 
 \subsection{Application of Kuznetsov's formula} 
 
 Now we apply Kuznetsov's formula followed by the spectral large sieve, as in section \ref{appkuz}. We provide brief details, only for \eqref{need11} and only for the resulting Maass form contribution $\mathcal{M}$ of Theorem \ref{kuznetsov}. We consider two cases.
 
 \subsubsection{Case I} Suppose that
\begin{align}
\label{eq2} \frac{\sqrt{N'M'q}}{Cp^2} \le 1. 
\end{align}
By Lemma \ref{transf-bounds} (a), for real $t_f$ we may restrict to $t_f\ll (qp)^\epsilon$, for which we have $\tilde{\phi}(t_f)\ll (qp)^\epsilon$. We need to prove the following analogue of \eqref{largesieveexpectation}:
\begin{align}
 \label{doneb4} q^{\vartheta+\epsilon} p^{\epsilon} \frac{NM}{C^3 q} \Big( \Big(1+\frac{N'}{p^4}\Big)N'\Big)^\half \Big( \Big(1+\frac{M'}{p^4}\Big)M'\Big)^\half  \ll (qp)^\epsilon p^2q^{-\half+\vartheta}.
\end{align}
We consider the cross terms on the left hand side one by one: We have
\begin{align}
 \label{doneb4c} q^{\vartheta+\epsilon} p^{\epsilon} \frac{NM}{C^3 q} (N')^\half (M')^\half = p^{2+\epsilon}q^{-\half+\vartheta+\epsilon} \Big(\frac{\sqrt{NM}}{C\sqrt{q}}\Big)^2 \frac{\sqrt{N'M'q}}{Cp^2} \ll p^{2+\epsilon}q^{-\half+\vartheta+\epsilon}
\end{align}
on using \eqref{ccase1} and \eqref{eq2}; we have
\begin{align}
\label{doneb4a} q^{\vartheta+\epsilon} p^{\epsilon} \frac{NM}{C^3 q} \Big( \frac{N'}{p^4}N'\Big)^\half\Big( \frac{M'}{p^4}M'\Big)^\half =  q^{-\half+\vartheta+\epsilon} p^{2+\epsilon}  \frac{ \sqrt{NM}}{C\sqrt{q}} \frac{ \sqrt{NM}}{qp^2} \Big( \frac{\sqrt{N'M'q}}{Cp^2}\Big)^2  \ll p^{2+\epsilon}q^{-\half+\vartheta+\epsilon}
 \end{align}
on using \eqref{ccase1}, \eqref{eq3} and \eqref{eq2}; and we have
\begin{align}
\label{doneb4b} q^{\vartheta+\epsilon} p^{\epsilon} \frac{NM}{C^3 q}  (N')^\half \Big( \frac{M'}{p^4}M'\Big)^\half=   q^{-\half+\vartheta+\epsilon}  p^{2+\epsilon}\frac{ \sqrt{NM}}{C\sqrt{q}}  \frac{\sqrt{N'M'q}}{Cp^2} \frac{\sqrt{M'M}}{Cp}  \frac{ \sqrt{N}}{p\sqrt{q} } \ll   q^{-\half+\vartheta+\epsilon}  p^{2+\epsilon}.
\end{align}
on using \eqref{ccase1},  \eqref{eq2}, \eqref{n'm'range}, and $N\le (qp)^\epsilon q p^2$.

If $t_f$ is imaginary, then $t_f\ll 1 $ and $\tilde{\phi}(t_f)\ll (qp)^\epsilon (\frac{\sqrt{N'M'q}}{Cp^2})^{-2\vartheta}$ by taking $Z=1$ in Lemma \ref{transf-bounds} (b). The large sieve for the contribution of the exceptional eigenvalues then works in the same way as above, except that on the right hand sides of \eqref{doneb4c}, \eqref{doneb4a}, and \eqref{doneb4b}, we replace a factor of $\frac{\sqrt{N'M'q}}{Cp^2}$ with $(\frac{\sqrt{N'M'q}}{Cp^2})^{1-2\vartheta}$. Since $1-2\vartheta>0$, we end up with the same final bound $O(q^{-\half+\vartheta+\epsilon}  p^{2+\epsilon})$.

\subsubsection{Case II} Suppose that
\[
\frac{\sqrt{N'M'q}}{Cp^2} > 1. 
\]
Let $T:=  \frac{\sqrt{N'M'q}}{Cp^2}$. By Lemma \ref{transf-bounds} (a), for real $t_f$ we may restrict to $t_f\ll (qp)^\epsilon T$, for which we have $\tilde{\phi}(t_f)\ll (qp)^\epsilon  T^{-1} $. We need to prove the following analogue of \eqref{largesieveexpectation}:
\[
q^{\vartheta+\epsilon} p^{\epsilon} \frac{NM}{C^3 q}  T^{-1}  \Big( \Big(T^2+\frac{N'}{p^4}\Big)N'\Big)^\half  \Big( \Big(T^2+\frac{M'}{p^4}\Big)M'\Big)^\half  \ll (qp)^\epsilon p^2q^{-\half+\vartheta}.
\]
We consider the cross terms on the left one by one: We have
\[
 q^{\vartheta+\epsilon} p^{\epsilon} \frac{NM}{C^3 q} T^{-1} (T^2 N')^\half (T^2 M')^\half  = p^{2+\epsilon}q^{-\half+\vartheta+\epsilon} \frac{N N'}{C^2p^2} \frac{M M'}{C^2p^2}   \ll p^{2+\epsilon}q^{-\half+\vartheta+\epsilon}
\]
on using \eqref{n'm'range}; we have
\[
 q^{\vartheta+\epsilon} p^{\epsilon} \frac{NM}{C^3 q} \Big( \frac{N'}{p^4}N'\Big)^\half\Big( \frac{M'}{p^4} M' \Big)^\half  \ll p^{2+\epsilon}q^{-\half+\vartheta+\epsilon}
 \]
as already seen in \eqref{doneb4a}; and we have
\[
 q^{\vartheta+\epsilon} p^{\epsilon} \frac{NM}{C^3 q} T^{-1} (T^2 N')^\half \Big( \frac{M'}{p^4}M'\Big)^\half \ll   q^{-\half+\vartheta+\epsilon}  p^{2+\epsilon}.
\]
as already seen in \eqref{doneb4b}.

If $t_f$ is imaginary, then $t_f\ll 1$ and $\tilde{\phi}(t_f)\ll (qp)^\epsilon$ by taking $Z=1$ in Lemma \ref{transf-bounds} (b). The large sieve for the contribution of the exceptional eigenvalues then amounts to proving \eqref{doneb4}, which we have already done.

\bibliographystyle{amsplain} 

\bibliography{twisted-fourth}

\end{document}